\documentclass[aos,nameyear,dvips]{arximspdf}
\usepackage{dcolumn}
\usepackage{graphicx}

\doi{10.1214/09-AOS752}
\volume{38}
\issue{4}
\pubyear{2010}
\firstpage{2118}
\lastpage{2144}

\makeatletter

\newcolumntype{d}[1]{D{.}{.}{#1}}

\newtheorem{theorem}{Theorem}
\newtheorem{lemma}{Lemma}

\newproclaim{remark}{Remark}

\makeatother

\begin{document}
\begin{frontmatter}

\title{Optimal rates of convergence for covariance matrix estimation}
\runtitle{Covariance matrix estimation}

\begin{aug}
\author[A]{\fnms{T. Tony} \snm{Cai}\corref{}\thanksref{t1}\ead[label=e1]{tcai@wharton.upenn.edu}},
\author[B]{\fnms{Cun-Hui} \snm{Zhang}\thanksref{t2}\ead[label=e2]{cunhui@stat.rutgers.edu}} and
\author[C]{\fnms{Harrison H.} \snm{Zhou}\thanksref{t3}\ead[label=e3]{huibin.zhou@yale.edu}}
\runauthor{T. T. Cai, C.-H. Zhang and H. H. Zhou}
\affiliation{University of Pennsylvania, Rutgers University and Yale
University}
\address[A]{T. T. Cai\\
Department of Statistics\\
The Wharton School\\
University of Pennsylvania\\
Philadelphia, Pennsylvania 19104-6302\\
USA\\
\printead{e1}}
\address[B]{C.-H. Zhang\\
Department of Statistics\\
504 Hill Center\\
Busch Campus\\
Rutgers University\\
Piscataway, New Jersey 08854-8019\\
USA\\
\printead{e2}}
\address[C]{H. H. Zhou\\
Department of Statistics\\
Yale University\\
P.O. Box 208290\\
New Haven, Connecticut 06520-8290\\
USA\\
\printead{e3}}
\end{aug}

\thankstext{t1}{Supported in part by NSF Grant DMS-06-04954.}

\thankstext{t2}{Supported in part by NSF Grants DMS-05-04387 and DMS-06-04571.}

\thankstext{t3}{Supported in part by NSF Career Award DMS-06-45676.}

\received{\smonth{10} \syear{2008}}
\revised{\smonth{9} \syear{2009}}

%
\begin{abstract}
Covariance matrix plays a central role in multivariate statistical analysis.
Significant advances have been made recently on developing both theory and
methodology for estimating large covariance matrices. However, a minimax
theory has yet been developed.
In this paper we establish the optimal rates of convergence for estimating
the covariance matrix under both the operator norm and Frobenius norm.
It is
shown that optimal procedures under the two norms are different and
consequently matrix estimation under the operator norm is fundamentally
different from vector estimation. The minimax upper bound is obtained by
constructing a special class of tapering estimators and by studying their
risk properties. A key step in obtaining the optimal rate of
convergence is
the derivation of the minimax lower bound. The technical analysis requires
new ideas that are quite different from those used in the more conventional
function/sequence estimation problems.
\end{abstract}

%

\end{frontmatter}

\section{Introduction}

Suppose we observe independent and identically distributed $p$-variate
random variables $\mathbf{X}_{1},\ldots,\mathbf{X}_{n}$ with covariance
matrix $\Sigma_{p\times p}$ and the goal is to estimate the unknown
matrix $%
\Sigma_{p\times p}$ based on the sample $ \{ \mathbf{X}_{i}\dvtx i =1,
\ldots, n \} $. This covariance matrix estimation problem is of
fundamental importance in multivariate analysis. A wide range of statistical
methodologies, including clustering analysis, principal component analysis,
linear and quadratic discriminant analysis, regression analysis,
require the
estimation of the covariance matrices. With \mbox{dramatic} advances in technology,
large high-dimensional data are now routinely collected in scientific
investigations. Examples include climate studies, gene expression arrays,
functional magnetic resonance imaging, risk management and portfolio
allocation and web search problems. In such settings, the standard and most
natural estimator, the sample covariance matrix, often performs poorly. See,
for example, Muirhead (\citeyear{M87}), Johnstone (\citeyear{J01}), Bickel and Levina
(\citeyear{BL08a}, \citeyear{BL08b})
and Fan, Fan and Lv (\citeyear{FFL08}).


Regularization methods, originally developed in nonparametric function
estimation, have recently been applied to estimate large covariance
matrices. These include banding method
in Wu and Pourahmadi (\citeyear{WP09}) and Bickel and Levina (\citeyear{BL08a}),
tapering in Furrer and Bengtsson (\citeyear{FB07}), thresholding in
Bickel and Levina (\citeyear{BL08b}) and El Karoui (\citeyear{K08}), penalized estimation in
Huang et al. (\citeyear{Huangetal06}), Lam and Fan (\citeyear{LF07}) and Rothman
et al. (\citeyear{Rothmanetal08}), regularizing principal components in
Johnstone and Lu (\citeyear{JL04}) and Zou, Hastie and Tibshirani (\citeyear{ZHT06}). Asymptotic
properties and convergence results have been given in several papers. In
particular, Bickel and Levina (\citeyear{BL08a}, \citeyear{BL08b}), El Karoui (\citeyear{K08}) and Lam and
Fan (\citeyear{LF07}) showed consistency of their estimators in operator norm and even
obtained explicit rates of convergence. However, it is not clear
whether any
of these rates of convergence are optimal.

Despite recent progress on covariance matrix estimation there has been
remarkably little fundamental theoretical study on optimal estimation.
In this paper, we establish the optimal rate of convergence for estimating
the covariance matrix as well as its inverse over a wide range of
classes of
covariance matrices. Both the operator norm and Frobenius norm are
considered. It is shown that optimal procedures for these two norms are
different and consequently matrix estimation under the operator norm is
fundamentally different from vector estimation. In addition, the results
also imply that the banding estimator given in Bickel and Levina
(\citeyear{BL08a}) is
sub-optimal under the operator norm and the performance can be significantly
improved.

We begin by considering optimal estimation of the covariance matrix
$\Sigma$
over a class of matrices that has been considered in Bickel and Levina
(\citeyear{BL08a}). Both minimax lower and upper bounds are derived. We write $%
a_{n}\asymp b_{n}$ if there are positive constants $c$ and $C$ independent
of $n$ such that $c\leq a_{n}/b_{n}\leq C$. For a matrix $A$ its operator
norm is defined as $ \Vert A \Vert=\sup_{ \Vert
x \Vert
_{2}=1} \Vert Ax \Vert_{2}$. We assume that $p\leq\exp
(
\gamma n ) $ for some constant $\gamma>0$. Combining the results given
in Section \ref{operatornorm.sec}, we have the following optimal rate of
convergence for estimating the covariance matrix under the operator norm.
\begin{theorem}
\label{MinimaxOpe} The minimax risk of estimating the covariance
matrix $%
\Sigma$ over the class $\mathcal{P}_{\alpha}$ given in (\ref{paraspace})
satisfies
%
%
\begin{equation}\label{rateOper}
\inf_{\hat{\Sigma}}\sup_{\mathcal{P}_{\alpha}}\mathbb{E}
\Vert\hat{%
\Sigma}-\Sigma\Vert^{2}\asymp\min\biggl\{ n^{-
{2\alpha}/({%
2\alpha+1})}+\frac{\log p}{n}, \frac{p}{n} \biggr\}.
\end{equation}
\end{theorem}

The minimax upper bound is obtained by constructing a class of tapering
estimators and by studying their risk properties. It is shown that the
estimator with the optimal choice of the tapering parameter attains the
optimal rate of convergence. In comparison to some existing methods in the
literature, the proposed procedure does not attempt to estimate each
row/column optimally as a vector. In fact, our procedure does not optimally
trade bias and variance for each row/column. As a vector estimator, it has
larger variance than squared bias for each row/column. In other words,
it is
undersmoothed as a vector.

A key step in obtaining the optimal rate of convergence is the
derivation of
the minimax lower bound. The lower bound is established by using a testing
argument, where at the core is a novel construction of a collection of least
favorable multivariate normal distributions and the application of Assouad's
lemma and Le Cam's method. The technical analysis requires ideas that are
quite different from those used in the more conventional function/sequence
estimation problems.

In addition to the asymptotic analysis, we also carry out a small simulation
study to investigate the finite sample performance of the proposed
estimator. The tapering estimator is easy to implement. The numerical
performance of the estimator is compared with that of the banding estimator
introduced in Bickel and Levina (\citeyear{BL08a}). The simulation study shows
that the
proposed estimator has good numerical performance; it nearly uniformly
outperforms the banding estimator.

The paper is organized as follows. In Section \ref{method.sec}, after basic
notation and definitions are introduced, we propose a tapering procedure
for the covariance matrix estimation. Section \ref{operatornorm.sec} derives
the optimal rate of convergence for estimation under the operator norm. The
upper bound is obtained by studying the properties of the tapering
estimators and the minimax lower bound is obtained by a testing argument.
Section \ref{Frobeniusnorm.sec} considers optimal estimation under the
Frobenius norm. The problem of estimating the inverse of a covariance matrix
is treated in Section \ref{inverse.sec}. Section \ref{simulation.sec}
investigates the numerical performance of our procedure by a simulation
study. The technical proofs of auxiliary lemmas are given in Section
\ref%
{sec.proofs}.

\section{Methodology}

\label{method.sec}

In this section we will introduce a tapering procedure for estimating the
covariance matrix $\Sigma_{p\times p}$ based on a random sample of $p$%
-variate observations $\mathbf{X}_{1},\ldots,\mathbf{X}_{n}$. The
properties of the tapering estimators under the operator norm and Frobenius
norm are then studied and used to establish the minimax upper bounds in
Sections \ref{operatornorm.sec} and \ref{Frobeniusnorm.sec}.

Given a random sample $\{\mathbf{X}_{1},\ldots,\mathbf{X}_{n}\}$
from a
population with covariance matrix $\Sigma=\Sigma_{p\times p}$, the sample
covariance matrix is
\[
\frac{1}{n-1}\sum_{l=1}^{n} ( \mathbf{X}_{l}-\mathbf{\bar
{X}} )
( \mathbf{X}_{l}-\mathbf{\bar{X}} ) ^{T},
\]
which is an unbiased estimate of $\Sigma$, and the maximum likelihood
estimator of $\Sigma$ is
%
%
\begin{equation} \label{MLE}
\Sigma^{\ast}=(\sigma_{ij}^{\ast})_{1\leq i,j\leq p}=\frac{1}{n}%
\sum_{l=1}^{n} ( \mathbf{X}_{l}-\mathbf{\bar{X}} )
( \mathbf{X}%
_{l}-\mathbf{\bar{X}} ) ^{T},
\end{equation}
when $\mathbf{X}_{l}$'s are normally distributed. These two estimators are
close to each other for large $n$. We shall construct estimators of the
covariance matrix $\Sigma$ by tapering the maximum likelihood
estimator $%
\Sigma^{\ast}$.

Following Bickel and Levina (\citeyear{BL08a}) we consider estimating the covariance
matrix $\Sigma_{p\times p}= ( \sigma_{ij} ) _{1\leq
i,j\leq p}$
over the following parameter space:
%
%
\begin{eqnarray}\label{paraspace}
\mathcal{F}_{\alpha}=\mathcal{F}_{\alpha} ( M_{0},M )
&=& \biggl\{
\Sigma\dvtx\max_{j}\sum_{i} \{ \vert\sigma_{ij} \vert
\dvtx \vert i-j \vert>k \} \leq Mk^{-\alpha}\nonumber\\[-8pt]\\[-8pt]
&&\hspace*{53.8pt}\mbox{for
all }k%
\mbox{, and }\lambda_{\max} ( \Sigma) \leq M_{0}
\biggr\},\nonumber
\end{eqnarray}
where $\lambda_{\max}(\Sigma)$ is the maximum eigenvalue of the
matrix $%
\Sigma$, and $\alpha>0$, $M>0$ and $M_{0}>0$. Note that the smallest
eigenvalue of any covariance matrix in the parameter space $F_{\alpha
}$ is
allowed to be $0$ which is more general than the assumption in (5)
of Bickel and Levina (\citeyear{BL08a}). The parameter $\alpha$ in (\ref{paraspace}),
which essentially specifies the rate of decay for the covariances
$\sigma
_{ij}$ as they move away from the diagonal, can be viewed as an analog of
the smoothness parameter in nonparametric function estimation problems. The
optimal rate of convergence for estimating $\Sigma$ over the parameter
space $\mathcal{F}_{\alpha} ( M_{0},M ) $ critically
depends on
the value of $\alpha$. Our estimators of the covariance matrix $\Sigma$
are constructed by tapering the maximum likelihood estimator (\ref
{MLE}) as
follows.

\subsection*{Estimation procedure}

\begin{figure}

\includegraphics{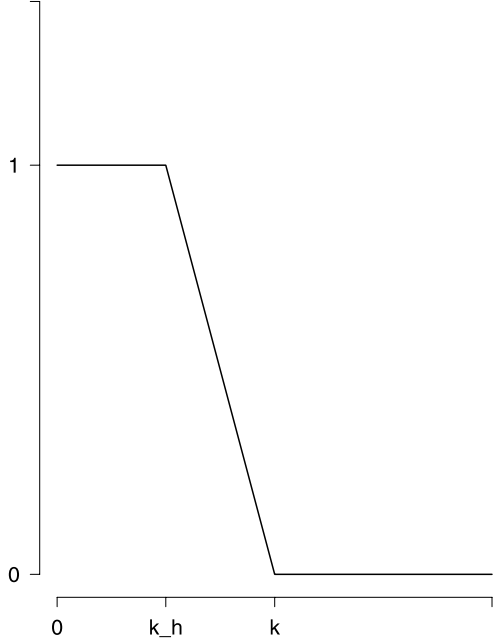}

\caption{The weights as a function of $|i-j|$.}
\label{weight.fig}
\end{figure}

For a given even integer $k$ with $1\leq k\leq p$, we define a tapering
estimator as
%
%
\begin{equation} \label{tapering.est}
\hat{\Sigma}=\hat{\Sigma}_{k}= ( w_{ij}\sigma_{ij}^{\ast
} )
_{p\times p},
\end{equation}
where $\sigma_{ij}^{\ast}$ are the entries in the maximum likelihood
estimator $\Sigma^{\ast}$ and the weights
%
%
\begin{equation} \label{wij}
w_{ij}=k_{h}^{-1}\{(k-|i-j|)_{+}-(k_{h}-|i-j|)_{+}\},
\end{equation}
where $k_{h}=k/2$. Without loss of generality we assume that $k$ is even.
Note that the weights $w_{ij}$ can be rewritten as
\[
w_{ij}= \cases{
1, &\quad when $|i-j|\leq k_{h}$, \cr
2-{\dfrac{|i-j|}{k_{h}}}, &\quad when $k_{h}<|i-j|<k$, \cr
0, &\quad otherwise.}
\]
See Figure \ref{weight.fig} for a plot of the weights $w_{ij}$ as a function
of $|i-j|$.

%

The tapering estimators are different from the banding estimators used in
Bickel and Levina (\citeyear{BL08a}). It is important to note that the tapering
estimator given in (\ref{tapering.est}) can be rewritten as a sum of many
small block matrices along the diagonal. This simple but important
observation is very useful for our technical arguments. Define the block
matrices
\[
M_{l}^{\ast( m ) }= ( \sigma_{ij}^{\ast}I \{
l\leq
i<l+m,l\leq j<l+m \} ) _{p\times p}
\]
and set
\[
S^{\ast( m ) }=\sum_{l=1-m}^{p}M_{l}^{\ast(
m ) }
\]
for all integers $1-m\leq l\leq p$ and $m\geq1$.
\begin{lemma}
\label{est} The tapering estimator $\hat{\Sigma}_{k}$ given in
(\ref{tapering.est}) can be written as
%
%
\begin{equation} \label{estimator}
\hat{\Sigma}_{k}=k_{h}^{-1} \bigl( S^{\ast( k )
}-S^{\ast(k_{h} ) } \bigr).
\end{equation}
\end{lemma}

It is clear that the performance of the estimator $\hat\Sigma_k$
depends on
the choice of the tapering parameter $k$. The optimal choice of $k$
critically depends on the norm under which the estimation error is measured.
We will study in the next two sections the rate of convergence of the
tapering estimator under both the operator norm and Frobenius norm.
Together with the minimax lower bounds derived in Sections \ref%
{operatornorm.sec} and \ref{Frobeniusnorm.sec}, the results show that a
tapering estimator with the optimal choice of $k$ attains the optimal rate
of convergence under these two norms.

\section{Rate optimality under the operator norm}
\label{operatornorm.sec}

In this section we will establish the optimal rate of convergence under the
operator norm. For $1\leq q\leq\infty$, the matrix $\ell_{q}$-norm
of a
matrix $A$ is defined by $\Vert A\Vert_{q}={\max_{\Vert x\Vert
_{q}=1}}\Vert
Ax\Vert_{q}$. The commonly used operator norm $\Vert\cdot\Vert$
coincides with the matrix $\ell_{2}$-norm $\Vert\cdot\Vert_{2}$.
For a
symmetric matrix $A$, it is known that the operator norm $ \Vert
A \Vert$ is equal to the largest magnitude of eigenvalues of $A$.
Hence it is also called the spectral norm. We will establish Theorem
\ref{MinimaxOpe} by deriving a minimax upper bound using the tapering estimator
and a matching minimax lower bound by a careful construction of a collection
of multivariate normal distributions and the application of Assouad's lemma
and Le Cam's method. We shall focus on the case $p\geq n^{{1}/({2\alpha+1%
})}$ in Sections~\ref{sec.uppbd} and~\ref{sec.lowbd}. The case of
$p<n^{{1}/({2\alpha+1})}$, which will be discussed in Section~\ref{sec.discussion},
is similar and slightly easier.

\subsection{Minimax upper bound under the operator norm}
\label{sec.uppbd}

We derive in this section the risk upper bound for the tapering estimators
defined in (\ref{estimator}) under the operator norm. Throughout the paper
we denote by $C$ a generic positive constant which may vary from place to
place but always depends only on indices $\alpha$, $M_{0}$ and $M$ of the
matrix family. We shall assume that the distribution of the $X_{i}$'s is
sub-Gaussian in the sense that there is $\rho>0$ such that
%
%
\begin{equation} \label{subGau}
\mathbb{P}\{|\mathbf{v}^T (\mathbf{X}_{1}-\mathbb{E}\mathbf{X}%
_{1})|>t\}\leq e^{-t^{2}\rho/2}\qquad\mbox{for all }t>0\mbox{ and }\Vert
\mathbf{%
v}\Vert_{2}=1.
\end{equation}
Let $\mathcal{P}_{\alpha}=\mathcal{P}_{\alpha} (
M_{0},M,\rho
) $ denote the set of distributions of $\mathbf{X}_{1}$ that
satisfy (%
\ref{paraspace}) and~(\ref{subGau}).
\begin{theorem}
\label{Operupper.bd.thm} The tapering estimator $\hat{\Sigma}_{k}$, defined
in (\ref{estimator}), of the covariance matrix $\Sigma_{p\times p}$
with $%
p\geq n^{{1}/({2\alpha+1})}$ satisfies
%
%
\begin{equation} \label{Operupbd.k}
\sup_{\mathcal{P}_{\alpha}}\mathbb{E} \Vert\hat{\Sigma
}_{k}-\Sigma
\Vert^{2}\leq C\frac{k+\log p}{n}+Ck^{-2\alpha}
\end{equation}
for $k=o ( n )$, $\log p=o ( n ) $ and some
constant $C>0$%
. In particular, the estimator $\hat{\Sigma}=\hat{\Sigma}_{k}$ with
$k=n^{{1}/({2\alpha+1})}$ satisfies
%
%
\begin{equation} \label{Operupbd}
\sup_{\mathcal{P}_{\alpha}}\mathbb{E} \Vert\hat{\Sigma
}-\Sigma
\Vert^{2}\leq Cn^{-{2\alpha}/({2\alpha+1})}+C\frac{\log p}{n}.
\end{equation}
\end{theorem}

From (\ref{Operupbd.k}) it is clear that the optimal choice of $k$ is of
order $n^{{1}/({2\alpha+1})}$. The upper bound given in (\ref{Operupbd})
is thus rate optimal among the class of the tapering estimators defined
in (%
\ref{estimator}). The minimax lower bound derived in Section \ref{sec.lowbd}
shows that the estimator $\hat{\Sigma}_{k}$ with $k=n^{{1}/({2\alpha+1})}$
is in fact rate optimal among all estimators.
\begin{pf*}{Proof of Theorem \protect\ref{Operupper.bd.thm}}
Note that $\Sigma^{\ast}$ is translation invariant and so is $\hat
{\Sigma}$. We
shall thus assume $\mu=0$ for the rest of the paper. Write%
\[
\Sigma^{\ast}=\frac{1}{n}\sum_{l=1}^{n} ( \mathbf
{X}_{l}-\mathbf{\bar{X%
}} ) ( \mathbf{X}_{l}-\mathbf{\bar{X}} )
^{T}=\frac{1}{n}%
\sum_{l=1}^{n}\mathbf{X}_{l}\mathbf{X}_{l}^{T}-\mathbf{\bar{X}\bar{X}}^{T},
\]
where $\mathbf{\bar{X}\bar{X}}^{T}$ is a higher order term (see
Remark \ref%
{higherorder} at the end of this section). In what follows we shall ignore
this negligible term and focus on the dominating term $\frac{1}{n}%
\sum_{l=1}^{n}\mathbf{X}_{l}\mathbf{X}_{l}^{T}$.

Set $\tilde{\Sigma}=\frac{1}{n}\sum_{l=1}^{n}\mathbf{X}_{l}\mathbf
{X}%
_{l}^{T} $ and write $\tilde{\Sigma}= ( \tilde{\sigma
}_{ij} )
_{1\leq i,j\leq p}$. Let
%
%
\begin{equation} \label{estimator1}
\breve{\Sigma}= ( \breve{\sigma}_{ij} ) _{1\leq i,j\leq
p}= (
w_{ij}\tilde{\sigma}_{ij} ) _{1\leq i,j\leq p}
\end{equation}
with $w_{ij}$ given in (\ref{wij}). Let $\mathbf{X}_{l}= (
X_{1}^{l},X_{2}^{l},\ldots,X_{p}^{l} ) ^{T}$. We then write
$\tilde{%
\sigma}_{ij}=\frac{1}{n}\sum_{l=1}^{n}X_{i}^{l}X_{j}^{l}$. It is
easy to see
%
%
\begin{eqnarray}
\label{meansc}
\mathbb{E}\tilde{\sigma}_{ij} &=&\sigma_{ij},\\
\label{varsc}
\hspace*{20pt}\operatorname{Var}(\tilde{\sigma}_{ij}) &=&\frac{1}{n}\operatorname{Var} (
X_{i}^{l}X_{j}^{l} ) \leq\frac{1}{n}\mathbb{E} (
X_{i}^{l}X_{j}^{l} ) ^{2}\leq\frac{1}{n}
[\mathbb E( X_{i}^{l} )^{4}]^{1/2}[\mathbb E ( X_{j}^{l} ) ^{4}]^{1/2}
\leq\frac{C}{n},\hspace*{-10pt}
\end{eqnarray}
that is, $\tilde{\sigma}_{ij}$ is an unbiased estimator of $\sigma
_{ij}$ with
a variance $O ( 1/n ) $.

We will first show that the variance part satisfies
%
%
\begin{equation} \label{varupp}
\mathbb{E} \Vert\breve{\Sigma}-\mathbb{E}\breve{\Sigma
} \Vert
^{2}\leq C\frac{k+\log p}{n}
\end{equation}
and the bias part satisfies
%
%
\begin{equation} \label{biasupp}
\Vert\mathbb{E}\breve{\Sigma}-\Sigma\Vert^{2}\leq
Ck^{-2\alpha}.
\end{equation}
It then follows immediately that
\[
\mathbb{E} \Vert\breve{\Sigma}-\Sigma\Vert^{2}\leq
2\mathbb{E}%
\Vert\breve{\Sigma}-\mathbb{E}\breve{\Sigma} \Vert
^{2}+2 \Vert\mathbb{E}\breve{\Sigma}-\Sigma\Vert
^{2}\leq
2C \biggl( \frac{k+\log p}{n}+k^{-2\alpha} \biggr).
\]
This proves (\ref{Operupbd.k}) and equation (\ref{Operupbd}) then
follows. Since $p\geq n^{{1}/({2\alpha+1})}$, we may choose
%
%
\begin{equation} \label{kope}
k=n^{{1}/({2\alpha+1})}
\end{equation}
and the estimator $\hat{\Sigma}$ with $k$ given in (\ref{kope}) satisfies
\[
\mathbb{E} \Vert\hat{\Sigma}-\Sigma\Vert^{2}\leq
2C \biggl( n^{-{2\alpha}/({2\alpha+1})}+\frac{\log p}{n} \biggr).
\]
Theorem \ref{Operupper.bd.thm} is then proved. 
\end{pf*}

We first prove the risk upper bound (\ref{biasupp}) for the bias part.
It is
well known that the operator norm of a symmetric matrix $A= (
a_{ij} ) _{p\times p}$ is bounded by its $\ell_{1}$ norm, that
is,%
\[
\Vert A \Vert\leq\Vert A\Vert_{1}=\max_{i=1,\ldots
,p}\sum_{j=1}^{p} \vert a_{ij} \vert
\]
[see, e.g., page 15 in Golub and Van Loan (\citeyear{GL83})]. This result was used in
Bickel and Levina (\citeyear{BL08a}, \citeyear{BL08b}) to obtain rates of convergence for their
proposed procedures under the operator norm (see discussions in Section
\ref%
{sec.discussion}). We bound the operator norm of the bias part $\mathbb
{E}%
\breve{\Sigma}-\Sigma$ by its $\ell_{1}$ norm. Since $\mathbb
{E}\tilde{%
\sigma}_{ij}=\sigma_{ij}$, we have%
\[
\mathbb{E}\breve{\Sigma}-\Sigma= \bigl( ( w_{ij}-1 )
\sigma
_{ij} \bigr) _{p\times p},
\]
where $w_{ij}\in[ 0,1 ] $ and is exactly $1$ when
$|i-j|\leq k$,
then%
\[
\Vert\mathbb{E}\breve{\Sigma}-\Sigma\Vert^{2}\leq
\biggl[
{\max_{i=1,\ldots,p}\sum_{j\dvtx|i-j|>k} }\vert\sigma_{ij}
\vert%
\biggr] ^{2}\leq M^{2}k^{-2\alpha}.
\]

Now we establish (\ref{varupp}) which is relatively complicated. The key
idea in the proof is to write the whole matrix as an average of matrices
which are sum of a large number of small disjoint block matrices, and for
each small block matrix the classical random matrix theory can be applied.
The following\vspace*{2pt} lemma shows that the operator norm of the random matrix $%
\breve{\Sigma}-\mathbb{E}\breve{\Sigma}$ is controlled by the
maximum of
operator norms of $p$ number of $k\times k$ random matrices. Let $%
M_{l}^{ ( m ) }= ( \tilde{\sigma}_{ij}I \{ l\leq
i<l+m,l\leq j<l+m \} ) _{p\times p}$. Define%
\[
N_{l}^{(m)}=\max_{1\leq l\leq p-m+1} \bigl\Vert M_{l}^{(m)}-\mathbb{E}
M_{l}^{(m)} \bigr\Vert.
\]
\begin{lemma}
\label{estbias}Let $\breve{\Sigma}$ be defined as in (\ref{estimator}). Then
\[
\Vert\breve{\Sigma}-\mathbb{E}\breve{\Sigma} \Vert
\leq
3N_{l}^{(m)}.
\]
\end{lemma}

For each small $m\times m$ random matrix with $m=k$, we control its operator
norm as follows.
\begin{lemma}
\label{estbiasbd}There is a constant $\rho_{1}>0$ such that
%
%
\begin{equation} \label{Mtail}
\mathbb{P} \bigl\{ N_{l}^{(m)}>x \bigr\} \leq2p5^{m}\exp(
-nx^{2}\rho
_{1} )
\end{equation}
for all $0<x<\rho_{1}$ and $1-m\leq l\leq p$.
\end{lemma}

With Lemmas \ref{estbias} and \ref{estbiasbd} we are now ready to
show the
variance bound (\ref{varupp}). By Lemma \ref{estbias} we have
\begin{eqnarray*}
\mathbb{E} \Vert\breve{\Sigma}-\mathbb{E}\breve{\Sigma
} \Vert^{2}
&\leq&9\mathbb{E} \bigl( N_{l}^{(m)} \bigr) ^{2}=9\mathbb{E} \bigl(
N_{l}^{(m)} \bigr) ^{2} \bigl[ I \bigl( N_{l}^{(m)}\leq x \bigr)
+I \bigl(
N_{l}^{(m)}>x \bigr) \bigr] \\
&\leq&9\bigl[ x^{2}+\mathbb{E} \bigl( N_{l}^{(m)} \bigr) ^{2}I \bigl(
N_{l}^{(m)}>x \bigr) \bigr].
\end{eqnarray*}
Note that $\Vert\mathbb{E}\breve{\Sigma}\Vert\leq\Vert
\Sigma\Vert$, which is bounded by a constant, and $\Vert\breve{%
\Sigma}\Vert\leq\Vert\breve{\Sigma}\Vert_{F}$. The
Cauchy--Schwarz inequality then implies%
\begin{eqnarray*}
\mathbb{E} \Vert\breve{\Sigma}-\mathbb{E}\breve{\Sigma
} \Vert^{2}
&\leq&C_{1} \bigl[ x^{2}+\mathbb{E} ( \Vert\breve{\Sigma
}%
\Vert_{F}^{2}+C ) I \bigl( N_{l}^{(m)}>x \bigr) \bigr]
\\
&\leq&C_{1} \bigl[ x^{2}+\sqrt{\mathbb{E} ( \Vert\breve
{\Sigma}%
\Vert_{F}+C ) ^{4}}\sqrt{\mathbb{P} \bigl(
N_{l}^{(m)}>x \bigr) }%
\bigr].
\end{eqnarray*}
Set $x=4\sqrt{\frac{\log p+m}{n\rho_{1}}}$. Then $x$ is bounded by
$\rho
_{1}$ as $n\rightarrow\infty$. From Lemma \ref{estbiasbd} we obtain%
%
%
\begin{eqnarray}\label{varupp1}
\mathbb{E} \Vert\breve{\Sigma}-\mathbb{E}\breve{\Sigma
} \Vert
^{2} &\leq& C \biggl[ \frac{\log p+m}{n}+p^{2}\cdot( p5^{m}\cdot
p^{-8}e^{-8m} ) ^{1/2} \biggr] \nonumber\\[-8pt]\\[-8pt]
&\leq& C_{1} \biggl( \frac{\log
p+m}{n}\biggr).\nonumber
\end{eqnarray}
\begin{remark}
\label{higherorder}
In the proof of Theorem \ref{Operupper.bd.thm}, the term
$\mathbf{\bar{X}\bar{X}}^{T}$ was ignored. It is not
difficult to see that this term has negligible contribution after tapering.
Let $H=\mathbf{\bar{X}\bar{X}}^{T}$ and $H= ( h_{ij} )
_{p\times p}$%
. Define
\[
H_{l}^{ ( m ) }= ( h_{ij}I \{ l\leq i<l+m,l\leq
j<l+m \} ) _{p\times p}.
\]
Similarly to Lemma \ref{estbiasbd}, it
can be shown that%
%
%
\begin{equation} \label{Htail}
\mathbb{P} \Bigl\{ \max_{1\leq l\leq p-m+1} \bigl\Vert
H_{l}^{(m)}-\mathbb{E}%
H_{l}^{(m)} \bigr\Vert>t \Bigr\} \leq2p5^{m}\exp( -nt\rho
_{2} )
\end{equation}
for all $0<t<\rho_{2}$ and $1-m\leq l\leq p$. Note that $\mathbb
{E}H=\frac{1%
}{n}\Sigma$, then%
\[
\mathbb{E} \Vert H \Vert^{2}\leq2\mathbb{E} \Vert
H-\mathbb{E}%
H \Vert^{2}+2 \Vert\mathbb{E}H \Vert^{2}\leq
2\mathbb{E}%
\Vert H-\mathbb{E}H \Vert^{2}+2M_{0}^{2}/n^{2}.
\]
Let $t=16\frac{\log p+m}{n\rho_{2}}$. From (\ref{Htail}) we have
\begin{eqnarray*}
\mathbb{E} \Vert H-\mathbb{E}H \Vert^{2} &\leq&
t^{2}+\mathbb{E}%
\Vert H-\mathbb{E}H \Vert^{2}I \Bigl( \max_{1\leq l\leq
p-m+1} \bigl\Vert H_{l}^{(m)}-\mathbb{E}H_{l}^{(m)} \bigr\Vert
>t \Bigr) \\
&=&t^{2}+o ( t^{2} ) \leq C \biggl( \frac{\log p+m}{n}
\biggr) ^{2}
\end{eqnarray*}
by similar arguments as for (\ref{varupp1}). Therefore $H$
has a
negligible contribution to the risk.
\end{remark}

\subsection{Lower bound under the operator norm}
\label{sec.lowbd}

Theorem \ref{Operupper.bd.thm} in Section \ref{sec.uppbd} shows that the
optimal tapering estimator attains the rate of convergence
$n^{-{2\alpha}/({2\alpha+1})}+\frac{\log p}{n}$. In this section we shall
show that
this rate of convergence is indeed optimal among all estimators by showing
that the upper bound in equation (\ref{Operupbd}) cannot be improved. More
specifically we shall show that the following minimax lower bound holds.
\begin{theorem}
\label{Operlower.bd.thm} Suppose $p\leq\exp( \gamma n )
$ for
some constant $\gamma>0$. The minimax risk for estimating the covariance
matrix $\Sigma$ over $\mathcal{P}_{\alpha}$ under the operator norm
satisfies
\[
\inf_{\hat{\Sigma}}\sup_{\mathcal{P}_{\alpha}}\mathbb{E}
\Vert\hat{%
\Sigma}-\Sigma\Vert^{2}\geq cn^{-{2\alpha}/({2\alpha+1})}+c\frac{\log p}{n}.
\]
\end{theorem}

The basic strategy underlying the proof of Theorem \ref
{Operlower.bd.thm} is
to carefully construct a finite collection of multivariate normal
distributions and calculate the total variation affinity between pairs of
probability measures in the collection.

We shall now define a parameter space that is appropriate for the minimax
lower bound argument. For given positive integers $k$ and $m$ with
$2k\leq p$
and $1\leq m\leq k$, define the $p\times p$ matrix
$B(m,k)=(b_{ij})_{p\times
p}$ with
\[
b_{ij}=I \{ i=m\mbox{ and }m+1\leq j\leq2k\mbox{, or }j=m\mbox
{ and }%
m+1\leq i\leq2k \}.
\]
Set $k=n^{{1}/({2\alpha+1})}$ and $a=k^{- ( \alpha+1 )
}$. We
then define the collection of $2^{k}$ covariance matrices as
%
%
\begin{equation} \label{f11}\quad
\mathcal{F}_{11}= \Biggl\{ \Sigma( \theta) \dvtx\Sigma
( \theta
) =I_{p}+\tau a\sum_{m=1}^{k}\theta_{m}B(m,k), \theta
=(\theta
_{m})\in\{ 0,1 \} ^{k} \Biggr\},
\end{equation}
where $I_{p}$ is the $p\times p$ identity matrix and $0<\tau
<2^{-\alpha
-1}M $. Without loss of generality we assume that $M_{0}>1$ and $\rho>1$.
Otherwise we replace $I_{p}$ in (\ref{f11}) by $\varepsilon I_{p}$
for $0<\varepsilon<\min\{ M_{0},\rho\} $. For $0<\tau
<2^{-\alpha-1}M$ it is easy to check that $\mathcal{F}_{11}\subset
\mathcal{%
F}_{\alpha}(M_{0},M)$ as $n\rightarrow\infty$. In addition to
$\mathcal{F}%
_{11}$ we also define a collection of diagonal matrices
%
%
\begin{equation} \label{f12}\qquad
\mathcal{F}_{12}= \Biggl\{ \Sigma_{m}\dvtx\Sigma_{m}=I_{p}+ \Biggl( \sqrt
{\frac{%
\tau}{n}\log p_{1}}I \{ i=j=m \} \Biggr) _{p\times p},
0\leq
m\leq p_{1} \Biggr\},
\end{equation}
where $p_{1}=\min\{ p,e^{n/2} \} $ and $0<\tau<\min
\{
( M_{0}-1 ) ^{2}, ( \rho-1 ) ^{2},1 \} $.
Let $%
\mathcal{F}_{1}=\mathcal{F}_{11}\cup\mathcal{F}_{12}$. It is clear
that $%
\mathcal{F}_{1}\subset\mathcal{F}_{\alpha}(M_{0},M)$.

We shall show below separately that the minimax risks over multivariate
normal distributions with covariance matrix in (\ref{f11}) and (\ref{f12})
satisfy
%
%
\begin{equation}\label{operlwbdA}
\inf_{\hat{\Sigma}}\sup_{\mathcal{F}_{11}}\mathbb{E} \Vert
\hat{\Sigma}%
-\Sigma\Vert^{2}\geq cn^{-{2\alpha}/({2\alpha+1})}
\end{equation}
and
%
%
\begin{equation} \label{operlwbdB}
\inf_{\hat{\Sigma}}\sup_{\mathcal{F}_{12}}\mathbb{E} \Vert
\hat{\Sigma}%
-\Sigma\Vert^{2}\geq c\frac{\log p}{n}
\end{equation}
for some constant $c>0$. Equations (\ref{operlwbdA}) and (\ref{operlwbdB})
together imply
%
%
\begin{equation} \label{operlwbd}
\inf_{\hat{\Sigma}}\sup_{\mathcal{F}_{1}}\mathbb{E} \Vert
\hat{\Sigma}%
-\Sigma\Vert^{2}\geq\frac{c}{2}\biggl(n^{-{2\alpha}/({2\alpha+1})}+%
\frac{\log p}{n}\biggr)
\end{equation}
for multivariate normal distributions and this proves Theorem \ref%
{Operlower.bd.thm}. We shall establish the lower bound (\ref
{operlwbdA}) by
using Assouad's lemma in Section \ref{sec.assouad} and the lower bound
(\ref%
{operlwbdB}) by using Le Cam's method and a two-point argument in
Section %
\ref{sec.lecom}.

\subsubsection{A lower bound by Assouad's lemma}
\label{sec.assouad}

The key technical tool to establish equation (\ref{operlwbdA}) is
Assouad's lemma in Assouad (\citeyear{A83}). It gives a lower bound for the maximum
risk over the parameter set $\Theta= \{ 0,1 \} ^{k}$ to the
problem of estimating an arbitrary quantity $\psi( \theta
) $,
belonging to a metric space with metric $d$. Let $H ( \theta
,\theta^{\prime} ) =\sum_{i=1}^{k} \vert\theta
_{i}-\theta
_{i}^{\prime} \vert$ be the Hamming distance on $ \{
0,1 \}
^{k}$, which counts the number of positions at which $\theta$ and $%
\theta^{\prime}$ differ. For two probability measures $P$ and $Q$ with
density $p$ and $q$ with respect to any common dominating measure $\mu$,
write the total variation affinity $ \Vert P\wedge Q \Vert
=\int
p\wedge q\,d\mu$. Assouad's lemma provides a minimax lower bound for
estimating $\psi(\theta)$.
\begin{lemma}[(Assouad)]
\label{Assouad} Let $\Theta= \{ 0,1 \} ^{k}$ and let $T$
be an
estimator based on an observation from a distribution in the collection
$%
\{ P_{\theta},\theta\in\Theta\} $. Then for all $s>0$
\[
\max_{\theta\in\Theta}2^{s}\mathbb{E}_{\theta}d^{s} ( T,\psi
(
\theta) ) \geq\min_{H ( \theta,\theta^{\prime
} )
\geq1}\frac{d^{s} ( \psi( \theta) ,\psi(
\theta
^{\prime} ) ) }{H ( \theta,\theta^{\prime}
)} \cdot
\frac{k}{2} \cdot{\min_{H ( \theta,\theta^{\prime} )
=1}} \Vert\mathbb{P}_{\theta}\wedge\mathbb{P}_{\theta^{\prime
}} \Vert.
\]
\end{lemma}

Assouad's lemma is connected to multiple comparisons. In total there are
$k$ comparisons. The lower bound has three factors. The first factor is
basically the minimum cost of making a mistake per comparison, and the last
factor is the lower bound for the total probability of making type I and
type II errors for each comparison, and $k/2$ is the expected number of
mistakes one makes when $\mathbb{P}_{\theta}$ and $\mathbb{P}%
_{\theta^{\prime}}$ are not distinguishable from each other when
$H (
\theta,\theta^{\prime} ) =1$.

We now prove the lower bound (\ref{operlwbdA}). Let $\mathbf
{X}_{1},\ldots,%
\mathbf{X}_{n}\stackrel{\mathrm{i.i.d.}}{\sim} N ( 0,\Sigma
( \theta)
) $ with $\Sigma( \theta) \in\mathcal{F}_{11}$. Denote
the joint distribution by $P_{\theta}$. Applying Assouad's lemma to the
parameter space $\mathcal{F}_{11}$, we have
%
%
\begin{eqnarray} \label{lwd.a}
&&\inf_{\hat{\Sigma}}\max_{\theta\in\{0,1\}^{k}}2^{2}E_{\theta
} \Vert
\hat{\Sigma}-\Sigma( \theta) \Vert^{2}\nonumber\\[-8pt]\\[-8pt]
&&\qquad\geq
\min_{H (
\theta,\theta^{\prime} ) \geq1}\frac{ \Vert\Sigma
( \theta
) -\Sigma( \theta^{\prime} ) \Vert
^{2}}{H (
\theta,\theta^{\prime} ) }\frac{k}{2}\min_{H ( \theta
,\theta^{\prime} ) =1} \Vert P_{\theta}\wedge
P_{\theta^{\prime}} \Vert.\nonumber
\end{eqnarray}
We shall state the bounds for the the first and third factors on the
right-hand side of (\ref{lwd.a}) in two lemmas. The proofs of these
lemmas are given in
Section \ref{sec.proofs}.
\begin{lemma}
\label{dffbd} Let $\Sigma( \theta) $ be defined as in
(\ref{f11}). Then for some constant $c>0$
\[
\min_{H ( \theta,\theta^{\prime} ) \geq1}\frac{
\Vert\Sigma
( \theta) -\Sigma( \theta^{\prime} )
\Vert
^{2}}{H ( \theta,\theta^{\prime} ) }\geq cka^{2}.
\]
\end{lemma}
\begin{lemma}
\label{affbd}Let $\mathbf{X}_{1},\ldots,\mathbf{X}_{n}\stackrel
{\mathit{i.i.d.}}{\sim}%
N ( 0,\Sigma( \theta) ) $ with $\Sigma
( \theta
) \in\mathcal{F}_{11}$. Denote the joint distribution by
$P_{\theta}$%
. Then for some constant $c>0$
\[
\min_{H ( \theta,\theta^{\prime} ) =1} \Vert
P_{\theta}\wedge
P_{\theta^{\prime}} \Vert\geq c.
\]
\end{lemma}

It then follows from Lemmas \ref{dffbd} and \ref{affbd} together,
with the
fact $k=n^{{1}/({2\alpha+1})}$,
\[
\max_{\Sigma(\theta)\in\mathcal{F}_{11}}2^{2}E_{\theta}
\Vert\hat{%
\Sigma}-\Sigma( \theta) \Vert^{2}\geq\frac
{c^{2}}{2}%
k^{2}a^{2}\geq c_{1}n^{-{2\alpha}/({2\alpha+1})}
\]
for some $c_{1}>0$.

\subsubsection{A lower bound using Le Cam's method}
\label{sec.lecom}

We now apply Le Cam's method to derive the lower bound (\ref
{operlwbdB}) for
the minimax risk. Let $X$ be an observation from a distribution in the
collection $ \{ P_{\theta},\theta\in\Theta\}$ where
$\Theta
= \{ \theta_{0},\theta_{1},\ldots,\theta_{p_{1}} \}$. Le Cam's
method, which is based on a two-point testing argument, gives a lower bound
for the maximum estimation risk over the parameter set $\Theta$. More
specifically, let $L$ be the loss function. Define $r ( \theta
_{0},\theta_{m} ) =\inf_{t} [ L ( t,\theta_{0} )
+L ( t,\theta_{m} ) ] $ and $r_{\min}=\inf_{1\leq
m\leq
p_{1}}r ( \theta_{0},\theta_{m} ) $, and denote $\bar
{\mathbb{P}}=%
\frac{1}{p_{1}}\sum_{m=1}^{p_{1}}\mathbb{P}_{\theta_{m}}$.
\begin{lemma}
\label{LeCam} Let $T$ be an estimator of $\theta$ based on an observation
from a distribution in the collection $ \{ P_{\theta},\theta\in
\Theta
= \{ \theta_{0},\theta_{1},\ldots,\theta_{p_{1}} \}
\} $,
then
\[
\sup_{\theta}\mathbb{E}L ( T,\theta) \geq\frac
{1}{2}r_{\min
} \Vert\mathbb{P}_{\theta_{0}}\wedge\bar{\mathbb{P}}
\Vert.
\]
\end{lemma}

We refer to Yu (\citeyear{Y97}) for more detailed discussions on Le Cam's method.

To apply Le Cam's method, we need to first construct a parameter set.
For $%
1\leq m\leq p_{1}$, let $\Sigma_{m}$ be a diagonal covariance matrix
with $%
\sigma_{mm}=1+\sqrt{\tau\frac{\log p_{1}}{n}}$, $\sigma_{ii}=1$
for $%
i\neq m$, and let $\Sigma_{0}$ be the identity matrix. Let $\mathbf{X}
_{l}= ( X_{1}^{l},X_{2}^{l},\ldots,X_{p}^{l} ) ^{T}\sim
N (
0,\Sigma_{m} ) $, and denote the joint density of $\mathbf{X}%
_{1},\ldots,\mathbf{X}_{n}$ by $f_{m}$, $0\leq m\leq p_{1}$ with $%
p_{1}=\max\{ p, e^{n/2} \} $, which can be written as
follows:%
\[
f_{m}=\prod_{1\leq i\leq n,1\leq j\leq p,j\neq m}\phi
_{1} (
x_{j}^{i} ) \cdot\prod_{1\leq i\leq n}\phi_{\sigma
_{mm}} ( x_{m}^{i} ),
\]
where $\phi_{\sigma}$, $\sigma=1$ or $\sigma_{mm}$, is the density
of $%
N ( 0,\sigma^{2} ) $. Denote by $f_0$ the joint density of $
\mathbf{X}_{1},\ldots,\mathbf{X}_{n}$ when $\mathbf{X}_{l} \sim
N (
0,\Sigma_{0} ) $.

Let $\theta_{m}=\Sigma_{m}$ for $0\leq m\leq p_{1}$ and the loss
function $%
L$ be the squared operator norm. It is easy to see $r ( \theta
_{0},\theta_{m} ) =\frac{1}{2}\tau\frac{\log p_{1}}{n}$ for
all $%
1\leq m\leq p_{1} $. Then the lower bound (\ref{operlwbdB}) follows
immediately from Lemma \ref{LeCam} if there is a constant $c>0 $ such that
%
%
\begin{equation} \label{affbd2}
\Vert\mathbb{P}_{\theta_{0}}\wedge\bar{\mathbb{P}}
\Vert\geq c.
\end{equation}

Note that for any two densities $q_{0}$ and $q_{1}$, $\int q_{0}\wedge
q_{1}\,d\mu=1-\frac{1}{2}\int\vert q_{0}-q_{1} \vert\, d\mu$,
and Jensen's inequality implies
\[
\biggl[ \int\vert q_{0}-q_{1} \vert \,d\mu\biggr]
^{2}= \biggl( \int
\biggl\vert\frac{q_{0}-q_{1}}{q_{1}} \biggr\vert q_{1}\,d\mu\biggr)
^{2}\leq
\int\frac{ ( q_{0}-q_{1} ) ^{2}}{q_{1}}\,d\mu=\int\frac
{q_{0}^{2}}{%
q_{1}}\,d\mu-1.
\]
Hence $\int q_{0}\wedge q_{1}\,d\mu\geq1-\frac{1}{2} ( \int\frac
{%
q_{0}^{2}}{q_{1}}\,d\mu-1 ) ^{1/2}$. To establish equation (\ref
{affbd2}%
), it thus suffices to show that $\int( \frac{1}{p_{1}}%
\sum_{m=1}^{p_{1}}f_{m} ) ^{2}/f_{0}\,d\mu-1\rightarrow0$, that
is,%
%
%
\begin{equation} \label{chisquarebd}
\frac{1}{p_{1}^{2}}\sum_{m=1}^{p_{1}}\int\frac
{f_{m}^{2}}{f_{0}}\,d\mu+\frac{%
1}{p_{1}^{2}}\sum_{m\neq j}\int\frac{f_{m}f_{j}}{f_{0}}\,d\mu
-1\rightarrow0.
\end{equation}
We now calculate $\int\frac{f_{m}f_{j}}{f_{0}}\,d\mu$. For $m\neq j$
it is
easy to see%
\[
\int\frac{f_{m}f_{j}}{f_{0}}\,d\mu-1=0.
\]
When $m=j$, we have%
\begin{eqnarray*}
\int\frac{f_{m}^{2}}{f_{0}}\,d\mu&=&\frac{ ( \sqrt{2\pi\sigma
_{mm}}%
) ^{-2n}}{ ( \sqrt{2\pi} ) ^{-n}}\prod
_{1\leq i\leq
n}\int\exp\biggl[ ( x_{m}^{i} ) ^{2} \biggl( -\frac
{1}{\sigma_{mm}}%
+\frac{1}{2} \biggr) \biggr] \,dx_{m}^{i} \\
&=& [ 1- ( 1-\sigma_{mm} ) ^{2} ] ^{-n/2}=
\biggl( 1-\tau
\frac{\log p_{1}}{n} \biggr) ^{-n/2}.
\end{eqnarray*}
Thus%
%
%
\begin{eqnarray} \label{affbd3}
&&\int\Biggl( \frac{1}{p_{1}}\sum_{m=1}^{p_{1}}f_{m} \Biggr)
^{2}\bigg/f_{0}\,d\mu-1\nonumber\\
&&\qquad=\frac{1}{p_{1}^{2}}\sum_{m=1}^{p_{1}} \biggl( \int\frac
{f_{m}^{2}}{f_{0}}%
\,d\mu-1 \biggr) \nonumber\\[-8pt]\\[-8pt]
&&\qquad\leq \frac{1}{p_{1}} \biggl( 1-\tau\frac{\log p_{1}}{n} \biggr)
^{-n/2}-\frac{1}{p_{1}} \nonumber\\
&&\qquad=\exp\biggl[ -\log p_{1}-\frac{n}{2}\log\biggl( 1-\tau\frac{\log
p_{1}}{n}%
\biggr) \biggr] -\frac{1}{p_{1}}\rightarrow0\nonumber
\end{eqnarray}
for $0<\tau<1$, where the last step follows from the inequality $\log
( 1-x ) \geq-2x$ for $0<x<1/2$. Equation (\ref{affbd3}), together
with Lemma \ref{LeCam}, now immediately implies the lower bound given in
(\ref{operlwbdB}). %
\begin{remark}
\label{logp}
In covariance matrix estimation literature, it is
commonly assumed that $\frac{\log p}{n}\rightarrow0$. See, for example,
Bickel and Levina (\citeyear{BL08a}). The lower bound given in this section implies
that this assumption is necessary for estimating the covariance matrix
consistently under the operator norm.
\end{remark}

\subsection{Discussion}
\label{sec.discussion}

Theorems \ref{Operupper.bd.thm} and \ref{Operlower.bd.thm} together show
that the minimax risk for estimating the covariance matrices over the
distribution space $\mathcal{P}_{\alpha}$ satisfies, for $p\geq
n^{{1}/({2\alpha+1})}$,
%
%
\begin{equation} \label{minimax.rate}
\inf_{\hat{\Sigma}}\sup_{\mathcal{P}_{\alpha}}\mathbb{E}
\Vert\hat{%
\Sigma}-\Sigma\Vert^{2}\asymp n^{-{2\alpha}/({2\alpha+1})}+\frac{\log p}{n}.
\end{equation}
The results also show that the tapering estimator $\hat{\Sigma}_{k}$ with
tapering parameter $k=n^{{1}/({2\alpha+1})}$ attains the optimal
rate of
convergence $n^{-{2\alpha}/({2\alpha+1})}+\frac{\log p}{n}$.

A few interesting points can be made on the optimal rate of convergence
$n^{-{2\alpha}/({2\alpha+1})}+\frac{\log p}{n}$. When the dimension
$p$ is
relatively small, that is, $\log p=o(n^{{1}/({2\alpha+1})})$, $p$
has no
effect on the convergence rate and the rate is purely driven by the
``smoothness'' parameter $\alpha$.
However, when $p$ is large, that is, $\log p\gg n^{{1}/({2\alpha+1})}$, $p$
plays a significant role in determining the minimax rate.

We should emphasize that the optimal choice of the tapering parameter $k
\asymp n^{{1}/({2\alpha+1})}$ is different from the optimal choice for
estimating the rows/columns as vectors under mean squared error loss.
Straightforward calculation shows that in the latter case the best
cutoff is
$k \asymp n^{{1}/({2(\alpha+1)})}$ so that the tradeoff between the
squared bias and the variance is optimal. With $k \asymp
n^{{1}/({2\alpha+1})}$, the tapering estimator has smaller squared bias than the
variance as
a vector estimator of each row/column.

It is also interesting to compare our results with those given in
Bickel and
Levina (\citeyear{BL08a}). A banding estimator with bandwidth $k = (\frac
{\log p}{n})^{{1}/({2(\alpha+1)})}$ was proposed and the rate of
convergence $%
(\frac{\log p}{n} )^{{\alpha}/({\alpha+1})}$ was
proved. It is
easy to see that the banding estimator given in Bickel and Levina
(\citeyear{BL08a}) is
not rate optimal. Take, for example, $\alpha=1/2$ and $p=e^{\sqrt
{n}}$. Their
rate is $n^{-{1/6}}$, while the optimal rate in Theorem \ref%
{MinimaxOpe} is $n^{-{1/2}}$.

It is instructive to take a closer look at the motivation behind the
construction of the banding estimator in Bickel and Levina (\citeyear{BL08a}). Let the
banding estimator be
%
%
\begin{equation} \label{bandestimator}
\hat{\Sigma}_{B}= ( \sigma_{ij}^{\ast}I \{ |i-j|\leq
k \})
\end{equation}
and denote $\hat{\Sigma}_{B}-\mathbb{E}\hat{\Sigma}_{B}$ by $V$,
and let $%
V= ( v_{ij} ) $. An important step in the proof of Theorem 1 in
Bickel and Levina (\citeyear{BL08a}) is to control the operator norm by the $\ell_{1}$
norm as follows:
\begin{eqnarray*}
\mathbb{E} \Vert\hat{\Sigma}_{B}-\mathbb{E}\hat{\Sigma
}_{B} \Vert
^{2} &\leq&\mathbb{E} \Vert\hat{\Sigma}_{B}-\mathbb{E}\hat
{\Sigma}%
_{B} \Vert_{1}^{2}=\mathbb{E} \biggl( \max_{j=1,\ldots
,p}\sum_{i} \vert v_{ij} \vert\biggr) ^{2} \\
&\leq&C \biggl( \frac{k}{\sqrt{n}}\sqrt{\log p} \biggr) ^{2}=C\frac
{k^{2}\log p%
}{n}.
\end{eqnarray*}
Note that $\mathbb{E} [ \vert v_{ij} \vert I \{
|i-j|\leq
k \} ] \asymp1/\sqrt{n}$, then $\mathbb{E}\sum
_{i} \vert
v_{ij} \vert\asymp k/\sqrt{n}$. It is then expected that
$\mathbb{E}%
( \max_{j=1,\ldots,p}\sum_{i} \vert v_{ij} \vert
)
^{2}\leq C ( \frac{k}{\sqrt{n}}\sqrt{\log p} ) ^{2}$ [see Bickel
and Levina (\citeyear{BL08a}) for details] and so%
\[
\mathbb{E} \Vert\breve{\Sigma}-\Sigma\Vert
_{1}^{2}\leq C\frac{%
k^{2}\log p}{n}+Ck^{-2\alpha}.
\]
An optimal tradeoff of $k$ is then $ ( \frac{\log p}{n} )
^{{1%
}/({2 ( \alpha+1 ) })}$ which implies a rate of $ ( \frac
{\log p}{%
n} ) ^{-{\alpha}/({\alpha+1})}$ in Theorem 1 in Bickel and Levina
(\citeyear{BL08a}). This rate is slower than the optimal rate
$n^{-{2\alpha}/({2\alpha+1})}+\frac{\log p}{n}$ in Theorem \ref{MinimaxOpe}.

We have considered the parameter space $\mathcal{F}_{\alpha}$ defined
in (%
\ref{paraspace}). Other similar parameter spaces can also be
considered. For example, in time series analysis it is often assumed
the covariance $%
\vert\sigma_{ij} \vert$ decays at the rate $ \vert
i-j \vert^{- ( \alpha+1 ) }$ for some $\alpha>0$. Consider
the collection of positive-definite symmetric matrices satisfying the
following conditions:
%
%
\begin{eqnarray} \label{paraspace.g}
\mathcal{G}_{\alpha}&=&\mathcal{G}_{\alpha}(M_{0},M_{1})\nonumber\\[-8pt]\\[-8pt]
&=& \bigl\{
\Sigma
\dvtx \vert\sigma_{ij} \vert\leq M_{1} \vert i-j
\vert
^{-(\alpha+1)}\mbox{ for }i\neq j \mbox{ and } \lambda_{\max
}(\Sigma
)\leq M_{0} \bigr\},\nonumber
\end{eqnarray}
where $\lambda_{\max}(\Sigma)$ is the maximum eigenvalues of the
matrix $%
\Sigma$. Note that $\mathcal{G}_{\alpha}(M_{0}$, $M_{1})$ is a subset
of $%
\mathcal{F}_{\alpha} ( M_{0},M ) $ as long as $M_{1}\leq
\alpha M$%
. Using virtually identical arguments one can show that
\[
\inf_{\hat{\Sigma}}\sup_{\mathcal{P}_{\alpha}^{\prime}}\mathbb
{E}%
\Vert\hat{\Sigma}-\Sigma\Vert^{2}\asymp n^{-{2\alpha}/({2\alpha+1})}+\frac{\log
p}{n}.
\]
Let $\mathcal{P}_{\alpha}^{\prime}=\mathcal{P}_{\alpha}^{\prime
} ( M_{0},M,\rho) $ denote the set of distributions of
$\mathbf{X}%
_{1}$ that satisfies (\ref{subGau}) and~(\ref{paraspace.g}).
\begin{remark}
\label{nonpd}
Both the tapering estimator proposed in this paper and
banding estimator given in Bickel and Levina (\citeyear{BL08a}) are not necessarily
positive-semidefinite. A practical proposal to avoid this would be to
project the estimator $\hat{\Sigma}$ to the space of
positive-semidefinite matrices under the operator norm. More specifically,
one may first diagonalize $\hat{\Sigma}$ and then replace negative
eigenvalues by~$0$. The resulting estimator is then positive-semidefinite.
\end{remark}

\subsubsection{The case of $p<n^{{1}/({2\alpha+1})}$}

We have focused on the case $p\geq n^{{1}/({2\alpha+1})}$ in
Sections \ref%
{sec.uppbd} and \ref{sec.lowbd}. The case of $p<n^{{1}/({2\alpha+1})}$
can be handled in a similar way. The main difference is that in this
case we
no longer have a tapering estimator $\hat{\Sigma}_{k}$ with
$k=n^{{1}/({2\alpha+1})}$ because $k>p$. Instead the maximum likelihood estimator $%
\Sigma^{\ast}$ can be used directly. It is easy to show in this case
%
%
\begin{equation} \label{Operupbd.mle}
\sup_{\mathcal{P}_{\alpha}}\mathbb{E} \Vert\Sigma^{\ast
}-\Sigma
\Vert^{2}\leq C\frac{p}{n}.
\end{equation}
The lower bound can also be obtained by the application of Assouad's lemma
and by using a parameter space that is similar to $\mathcal{F}_{11}$.
To be
more specific, for an integer $1\leq m\leq p/2$, define the $p\times p$
matrix $B_{m}=(b_{ij})_{p\times p}$ with
\[
b_{ij}=I \{ i=m\mbox{ and }m+1\leq j\leq p\mbox{, or }j=m\mbox{
and }%
m+1\leq i\leq p \}.
\]
Define the collection of $2^{p/2}$ covariance matrices as
%
%
\begin{equation} \label{fstar}
\mathcal{F}^{\ast}= \Biggl\{ \Sigma( \theta) \dvtx\Sigma
(
\theta) =I_{p}+\tau\frac{1}{\sqrt{np}}\sum
_{m=1}^{p/2}\theta
_{m}B(m,k), \theta=(\theta_{m})\in\{ 0,1 \}
^{p/2} \Biggr\}.\hspace*{-28pt}
\end{equation}
Since $p<n^{{1}/({2\alpha+1})}$, then $\frac{1}{\sqrt
{np}}<2^{\alpha
+1/2}p^{- ( \alpha+1 ) }$. Again it is easy to check
$\mathcal{F}%
^{\ast}\subset\mathcal{F}_{\alpha}(M_{0},M)$ when $0<\tau
<2^{-\alpha
-1}M $. The following lower bound then follows from the same argument
as in
Section \ref{sec.assouad}:
%
%
\begin{equation} \label{operlwbd.mle}
\inf_{\hat{\Sigma}}\sup_{\mathcal{F}^{\ast}}\mathbb{E}
\Vert\hat{\Sigma%
}-\Sigma\Vert^{2}\geq cp \biggl( \frac{1}{\sqrt{np}} \biggr)
^{2}\cdot
\frac{p}{2}\cdot c_{1}\geq c_{2}\frac{p}{n}.
\end{equation}
Equations (\ref{Operupbd.mle}) and (\ref{operlwbd.mle}) together
yield the
minimax rate of convergence for the case $p\leq n^{{1}/({2\alpha+1})}$,
%
%
\begin{equation} \label{operbd.mle}
\inf_{\hat{\Sigma}}\sup_{\mathcal{P}_{\alpha}}\mathbb{E}
\Vert\hat{%
\Sigma}-\Sigma\Vert^{2}\asymp\frac{p}{n}.
\end{equation}
This, together with equation (\ref{minimax.rate}), gives the optimal
rate of
convergence:
%
%
\begin{equation}
\inf_{\hat{\Sigma}}\sup_{\mathcal{P}_{\alpha}}\mathbb{E}
\Vert\hat{%
\Sigma}-\Sigma\Vert^{2}\asymp\min\biggl\{ n^{-{2\alpha}/({2\alpha+1})}
+\frac{\log p}{n}, \frac{p}{n} \biggr\}.
\end{equation}

\section{Rate optimality under the Frobenius norm}
\label{Frobeniusnorm.sec}

In addition to the operator norm, the Frobenius norm is another commonly
used matrix norm. The Frobenius norm is used in defining the numerical rank
of a matrix which is useful in many applications, such as the principle
component analysis. See, for example, Rudelson and Vershynin (\citeyear{RV07}). The
Frobenius norm has also been used in the literature for measuring the
accuracy of a covariance matrix estimator. See, for example, Lam and
Fan (\citeyear{LF07}) and
Ravikumar et al. (\citeyear{Ravikumaretal08}). In this section we consider the optimal rate of
convergence for covariance matrix estimation under the Frobenius norm. The
Frobenius norm of a matrix $A= ( a_{ij} ) $ is defined as
the $\ell
_{2}$ vector norm of all entries in the matrix
\[
\Vert A \Vert_{F}=\sqrt{\sum_{i,j}a_{ij}^{2}}.
\]
This is equivalent to treating the matrix $A$ as a vector of length $p^{2}$.
It is easy to see that the operator norm is bounded by the Frobenius norm,
that is, $ \Vert A \Vert\leq\Vert A \Vert_{F}$.

The following theorem gives the minimax rate of convergence for estimating
the covariance matrix $\Sigma$ under the Frobenius norm based on the sample
$ \{ \mathbf{X}_{1},\ldots,\mathbf{X}_{n} \} $.
\begin{theorem}
\label{MinimaxFro}The minimax risk under the Frobenius norm satisfies
%
%
\begin{eqnarray}
\inf_{\hat{\Sigma}}\sup_{\mathcal{P}_{\alpha}}\mathbb{E}\frac
{1}{p}%
\Vert\hat{\Sigma}-\Sigma\Vert_{F}^{2}&\asymp&\inf_{\hat{\Sigma}}\sup_{\mathcal{P}^{\prime}_{\alpha}}\mathbb
{E}\frac{1}{p}%
\Vert\hat{\Sigma}-\Sigma\Vert_{F}^{2}\nonumber\\[-8pt]\\[-8pt]
&\asymp&\min
\biggl\{ n^{-%
({2\alpha+1})/({2 ( \alpha+1 ) })}, \frac{p}{n} \biggr\}.\nonumber
\end{eqnarray}
%
%
\end{theorem}

We shall establish below separately the minimax upper bound and minimax
lower bound.

\subsection{Upper bound under the Frobenius norm}

We will only prove the upper bound for the distribution set $\mathcal
{P}%
_{\alpha}^{\prime}$ given in (\ref{paraspace.g}). The
proof for
the parameter space $\mathcal{P}_{\alpha}$ is slightly more involved by thresholding procedures
as in Wavelet estimation. The minimax upper
bound is derived by again considering the tapering estimator (\ref%
{tapering.est}). Under the Frobenius norm the risk function is separable.
The risk of the tapering estimator can be bounded separately under the
squared $\ell_{2}$ loss for each row/column. This method has been commonly
used in nonparametric function estimation using orthogonal basis expansions.
Since
\[
\mathbb{E}\tilde{\sigma}_{ij}=\sigma_{ij}\quad\mbox{and}\quad
\operatorname{Var}(\tilde{%
\sigma}_{ij})\leq\frac{C}{n}
\]
for the tapering estimator (\ref{tapering.est}), we have
\[
\mathbb{E} ( w_{ij}\tilde{\sigma}_{ij}-\sigma_{ij} )
^{2}\leq
( 1-w_{ij} ) ^{2}\sigma_{ij}^{2}+w_{ij}^{2}\frac{C}{n}.
\]
It can be seen easily that%
\begin{eqnarray*}
\frac{1}{p}\mathbb{E} \Vert\breve{\Sigma}-\Sigma\Vert
_{F}^{2} &\leq& \frac{1}{p}\sum_{ \{ ( i,j )
\dvtx k_{h}<|i-j| \}
}\sigma_{ij}^{2}+\frac{1}{p}\sum_{ \{ ( i,j )
\dvtx|i-j|\leq
k \} } \biggl[ ( 1-w_{ij} ) ^{2}\sigma
_{ij}^{2}+w_{ij}^{2}%
\frac{C}{n} \biggr] \\
&\equiv& R_{1}+R_{2}.
\end{eqnarray*}
The assumption $\lambda_{\max} ( \Sigma) \leq M_{0}$ implies
that $\sigma_{ii}\leq M_{0}$ for all $i$. Since $ \vert\sigma
_{ij} \vert$ is also uniformly bounded for all $i\neq j$ from
assumption (\ref{paraspace.g}), we immediately have $R_{2}\leq C\frac{k}{n}$.

It is easy to show that%
%
%
\begin{equation} \label{Frouppbd1}
\frac{1}{p}\sum_{ \{ ( i,j ) \dvtx k<|i-j| \}
}\sigma
_{ij}^{2}\leq Ck^{-2\alpha-1},
\end{equation}
where $ \vert\sigma_{ij} \vert\leq C_{1} \vert
i-j \vert^{- ( \alpha+1 ) }$ for all $i\neq j$. Thus%
%
%
\begin{equation} \label{rateFro}
\mathbb{E}\frac{1}{p} \Vert\breve{\Sigma}-\Sigma\Vert
_{F}^{2}\leq Ck^{-2\alpha-1}+C\frac{k}{n}\leq C_{2}n^{-({2\alpha
+1})/({%
2 ( \alpha+1 ) })}
\end{equation}
by choosing
%
%
\begin{equation} \label{kFro}
k=n^{{1}/({2 ( \alpha+1 ) })}
\end{equation}
if $n^{{1}/({2 ( \alpha+1 ) })}\leq p$, which is
different from
the choice of $k$ for the operator norm in~(\ref{kope}). If
$n^{%
{1}/({2 ( \alpha+1 ) })}>p$, we will choose $k=p$, then
the bias
part is $0$ and consequently
\[
\mathbb{E}\frac{1}{p} \Vert\breve{\Sigma}-\Sigma\Vert
_{F}^{2}\leq C\frac{p}{n}.%
\]
\begin{remark}
For the parameter space $\mathcal{P}^{\prime}_{\alpha}$,\vspace*{2pt}
under the Frobenius norm the optimal tapering parameter $k$ is of
the order $n^{{1}/({2(\alpha+1)})}$. The rate of convergence
of the tapering estimator with $k\asymp n^{{1}/({2(\alpha+1)})}$
under the operator norm is
\[
\frac{\log p}{n}+n^{-{\alpha}/({\alpha+1})},
\]
which is slower than $n^{-{2\alpha}/({2\alpha+1})}+\frac{\log
p}{n}$ in (\ref{rateOper}). Similarly, the optimal procedure under the
operator norm is not rate optimal under the Frobenius norm. Therefore, the
optimal choice of the tapering parameter $k$ critically depends on the norm
under which the estimation accuracy is measured.
\end{remark}
\begin{remark}
Similarly for $\mathcal{P}^{\prime}_{\alpha}$, it can be shown that under the Frobenius norm the banding
estimator with $k\asymp n^{{1}/({2(\alpha+1)})}$ is rate optimal. Under
the operator norm, Bickel and Levina (\citeyear{BL08a}) chose $k\asymp(
\frac{%
\log p}{n} ) ^{{1}/({2(\alpha+1)})}$ for the banding
estimator which is close to $n^{{1}/({2(\alpha+1)})}$
up to a logarithmic factor of $p$. On the other hand, it can be shown
that for the parameter space $\mathcal{P}_{\alpha}$ no linear estimator
can achieve the optimal convergence rate under the Frobenius norm.
\end{remark}
%

\subsection{Lower bound under the Frobenius norm}

It is sufficient to establish the lower bound for the parameter space
$\mathcal{P}^{\prime}_{\alpha}$ given in (\ref{paraspace.g}).
Again the argument for $\mathcal{P}_{\alpha}$ is similar. As in the
case of
estimation under the operator norm, we need to construct a finite collection
of multivariate normal distributions with a parameter space $\mathcal
{G}%
_{2}\subset\mathcal{G}_{\alpha}$ such that%
\[
\inf_{\hat{\Sigma}}\sup_{\mathcal{G}_{2}}\mathbb{E}\frac
{1}{p} \Vert
\hat{\Sigma}-\Sigma\Vert_{F}^{2}\geq c\frac{k}{n}
\]
for some $c>0$ when $k=\min\{ n^{{1}/({2(\alpha+1)})}, p/2 \} $.

We construct $\mathcal{G}_{2}$ as follows. Let $0<\tau<M$ be a constant.
Define
\begin{eqnarray*}
&&\mathcal{G}_{2}= \bigl\{ \Sigma( \theta) \dvtx \Sigma
( \theta
) =I+ \bigl( \theta_{ij}\tau n^{-{{1/2}}}I \{ 1\leq
\vert i-j \vert\leq k \} \bigr) _{p\times p},\\
&&\hspace*{172.7pt}\mbox{for }%
\theta_{ij}=\theta_{ji}=0\mbox{ or }1 \bigr\}.
\end{eqnarray*}
It is easy to verify that $\mathcal{G}_{2}\subset\mathcal{G}_{\alpha
}$ as $%
n\rightarrow\infty$. Note that $\theta\mathbf{\in}\Theta= \{
0,1 \} ^{kp-k ( k+1 ) /2}$.

Applying Assouad's lemma with $d$ the Frobenius norm and $s=2$ to the
parameter space $\mathcal{G}_{2}$, we have
\begin{eqnarray*}
&&\max_{\theta\in\mathcal{G}_{2}}2^{2}E_{\theta}\frac{1}{p}
\Vert\hat{%
\Sigma}-\Sigma( \theta) \Vert_{F}^{2}\\
&&\qquad\geq\min
_{H (
\theta,\theta^{\prime} ) \geq1}{\frac{{1/p} \Vert
\Sigma
( \theta) -\Sigma( \theta^{\prime} )
\Vert
_{F}^{2}}{H ( \theta,\theta^{\prime} ) }\frac{kp-k (
k+1 ) /2}{2}\min_{H ( \theta,\theta^{\prime} )
=1}} \Vert P_{\theta}\wedge P_{\theta^{\prime}} \Vert.
\end{eqnarray*}
Note that%
\begin{eqnarray*}
\min_{H ( \theta,\theta^{\prime} ) \geq1}\frac
{1}{p}\frac{%
\Vert\Sigma( \theta) -\Sigma( \theta
^{\prime
} ) \Vert_{F}^{2}}{H ( \theta,\theta^{\prime
} ) }%
&=&\min_{H ( \theta,\theta^{\prime} ) \geq1}\frac
{1}{p}\frac{%
[ \tau n^{-{{1/2}}} ] ^{2}\sum\vert\theta
_{ij}-\theta_{ij}^{\prime} \vert^{2}}{H ( \theta,\theta
^{\prime} ) }\\
&=&\frac{\tau^{2}}{p}n^{-1}.
\end{eqnarray*}
It is easy to see that%
\[
\frac{kp-k ( k+1 ) /2}{2}\asymp kp.
\]
\begin{lemma}
\label{affbd1}Let $P_{\theta}$ be the joint distribution of $\mathbf
{X}%
_{1},\ldots,\mathbf{X}_{n}\stackrel{\mathit{i.i.d.}}{\sim}N (
0,\Sigma(
\theta) ) $ with $\Sigma( \theta) \in
\mathcal{G}%
_{2}$. Then for some constant $c_{1}>0$ we have%
\[
{\min_{H ( \theta,\theta^{\prime} ) =1} }\Vert
P_{\theta
}\wedge P_{\theta^{\prime}} \Vert\geq c_{1}.
\]
\end{lemma}

We omit the proof of this lemma. It is very similar to and simpler than the
proof of Lemma \ref{affbd}.

From Lemma \ref{affbd1} we have for some $c>0$%
%
%
\begin{equation} \label{Testtwo}
{\min_{H ( \theta,\theta^{\prime} ) =1}} \Vert
P_{\theta
}\wedge P_{\theta^{\prime}} \Vert\geq c
\end{equation}
thus%
\[
\max_{\theta\in\mathcal{G}_{2}}2^{2}E_{\theta}\frac{1}{p}
\Vert\hat{%
\Sigma}-\Sigma( \theta) \Vert_{F}^{2}\geq c\min
\biggl\{
n^{-({2\alpha+1})/({2 ( \alpha+1 ) })}, \frac
{p}{n} \biggr\},
\]
which implies that the rate obtained in (\ref{rateFro}) is optimal.

\section{Estimation of the inverse covariance matrix}
\label{inverse.sec}

The inverse of the covariance matrix $\Sigma^{-1}$ is of significant
interest in many statistical applications. The results and analysis
given in
Section \ref{operatornorm.sec} can be used to derive the optimal rate of
convergence for estimating $\Sigma^{-1}$ under the operator norm.

For estimating the inverse covariance matrix $\Sigma^{-1}$ we require the
minimum eigenvalue of $\Sigma$ to be bounded away from zero. For
$\delta>0$%
, we define%
%
%
\begin{equation} \label{eiglw}
L_{\delta}= \{ \Sigma\dvtx\lambda_{\min} ( \Sigma)
\geq
\delta\}.
\end{equation}
Let $\mathcal{\tilde{P}}_{\alpha}=\mathcal{\tilde{P}}_{\alpha
} (
M_{0},M,\rho,\delta) $ denote the set of distributions of
$\mathbf{X}%
_{1}$ that satisfy (\ref{paraspace}), (\ref{subGau}) and (\ref
{eiglw}), and
similarly, distributions in $\mathcal{\tilde{P}}^{\prime}_{\alpha}=
\mathcal{\tilde{P}}^{\prime}_{\alpha} (M_{0},M,\rho,\delta
) $
satisfy (\ref{subGau}), (\ref{paraspace.g}) and (\ref{eiglw}).

The following theorem gives the minimax rate of convergence for
estimating $%
\Sigma^{-1}$.
\begin{theorem}
\label{MinimaxInverseOpe} The minimax risk of estimating the inverse
covariance matrix $\Sigma^{-1}$ satisfies
%
%
\begin{equation} \label{rateInverseOper}
\inf_{\hat{\Sigma}}\sup_{\mathcal{\tilde{P}}}\mathbb{E}
\Vert\hat{\Sigma%
}^{-1}-\Sigma^{-1} \Vert^{2}\asymp\min\biggl\{ n^{-{2\alpha}/({2\alpha+1})}
+\frac{\log p}{n}, \frac{p}{n} \biggr\},
\end{equation}
where $\mathcal{\tilde{P}}$ denotes either $\mathcal{\tilde
{P}}_{\alpha}$
or $\mathcal{\tilde{P}}^{\prime}_{\alpha}$.
\end{theorem}
\begin{pf}
We shall focus on the case $p\geq n^{{1}/({2\alpha+1})}$.
The proof for the case of $p<n^{{1}/({2\alpha+1})}$ is
similar. To establish the upper bound, note that
\[
\hat{\Sigma}^{-1}-\Sigma^{-1}=\hat{\Sigma}^{-1} ( \Sigma
-\hat{\Sigma}%
) \Sigma^{-1},
\]
then%
\[
\Vert\hat{\Sigma}^{-1}-\Sigma^{-1} \Vert^{2}=
\Vert\hat{%
\Sigma}^{-1} ( \Sigma-\hat{\Sigma} ) \Sigma^{-1}
\Vert
^{2}\leq\Vert\hat{\Sigma}^{-1} \Vert^{2} \Vert
\Sigma-\hat{%
\Sigma} \Vert^{2} \Vert\Sigma^{-1} \Vert^{2}.
\]
It follows from assumption (\ref{paraspace}) that $ \Vert\Sigma
^{-1} \Vert^{2}\leq\delta^{-2}$. Note that
$\mathbb{P}%
\{ \Vert\breve{\Sigma}-\mathbb{E}\breve{\Sigma}
\Vert
^{2}>\epsilon\} \leq4p5^{m}\exp( -n\epsilon^{2}\rho
_{1} ) $ for any $\epsilon>0$ which decays faster than any polynomial
of $n$ as shown in the proof of Lemmas \ref{estbias} and \ref{estbiasbd}.
Let $\lambda_{\min} ( \breve{\Sigma} ) $ and $\lambda
_{\min
} ( \mathbb{E}\breve{\Sigma} ) $ be the smallest
eigenvalues of $%
\breve{\Sigma}$ and $\mathbb{E}\breve{\Sigma}$, respectively.
Then\break
$\mathbb{P}%
( \lambda_{\min} ( \breve{\Sigma} ) \leq\lambda
_{\min
} ( \mathbb{E}\breve{\Sigma} ) -\epsilon^{1/2} )
\geq\mathbb{%
P} ( \vert\lambda_{\min} ( \breve{\Sigma} )
-\lambda
_{\min} ( \mathbb{E}\breve{\Sigma} ) \vert\geq
\epsilon
^{1/2} ) $ decays faster than any polynomial of $n$. Let
$0<\epsilon<%
[ \lambda_{\min} ( \mathbb{E}\breve{\Sigma} )
/2 ] ^{2}$
and $c=1/ [ \lambda_{\min} ( \mathbb{E}\breve{\Sigma
} )
-\epsilon^{1/2} ] $, then $\mathbb{P} ( \Vert\hat
{\Sigma}%
^{-1} \Vert\geq c ) $ decays faster than any polynomial of $n$.
Therefore,
\begin{eqnarray*}
\mathbb{E} \Vert\hat{\Sigma}^{-1}-\Sigma^{-1} \Vert
^{2} &\leq&
\biggl( \frac{c}{\delta} \biggr) ^{2}\mathbb{E} \Vert\Sigma
-\hat{\Sigma%
} \Vert^{2}\\
&&{} + \mathbb{E} [ \Vert\hat{\Sigma
}^{-1} \Vert
^{2} \Vert\Sigma-\hat{\Sigma} \Vert^{2} \Vert
\Sigma
^{-1} \Vert^{2}I ( \Vert\hat{\Sigma}^{-1}
\Vert\geq
c ) ] \\
&\leq& C\min\biggl\{ n^{-{2\alpha}/({2\alpha+1})}+\frac{\log
p}{n}, %
\frac{p}{n} \biggr\}.
\end{eqnarray*}

The proof of the lower bound is almost identical to that of Theorem
\ref%
{MinimaxOpe} except that here we need to show%
\[
\min_{H ( \theta,\theta^{\prime} ) \geq1}\frac{
\Vert
\Sigma^{-1} ( \theta) -\Sigma^{-1} ( \theta
^{\prime
} ) \Vert^{2}}{H ( \theta,\theta^{\prime} )
}\geq
cka^{2}
\]
instead of Lemma \ref{dffbd}. For a positive definite matrix $A$, let $
\lambda_{\min} ( A ) $ denote the minimum eigenvalue of
$A$. Since%
\[
\Sigma^{-1} ( \theta) -\Sigma^{-1} ( \theta
^{\prime
} ) =\Sigma^{-1} ( \theta^{\prime} ) \bigl( \Sigma
(
\theta) -\Sigma( \theta^{\prime} ) \bigr)
\Sigma
^{-1} ( \theta) ,
\]
we have
\[
\Vert\Sigma^{-1} ( \theta) -\Sigma^{-1} (
\theta
^{\prime} ) \Vert\geq\lambda_{\min} ( \Sigma
^{-1} (
\theta) ) \lambda_{\min} ( \Sigma^{-1} (
\theta
^{\prime} ) ) \Vert\Sigma( \theta)
-\Sigma
( \theta^{\prime} ) \Vert.
\]
Note that
\[
\lambda_{\min} ( \Sigma^{-1} ( \theta) ) >1/M_{0},\qquad
\lambda_{\min} ( \Sigma^{-1} ( \theta^{\prime} )
) >1/M_{0},
\]
then Lemma \ref{dffbd} implies
\[
\min_{H ( \theta,\theta^{\prime} ) \geq1}\frac{
\Vert
\Sigma^{-1} ( \theta) -\Sigma^{-1} ( \theta
^{\prime
} ) \Vert^{2}}{H ( \theta,\theta^{\prime} )
}\geq
M_{0}^{-4}\min_{H ( \theta,\theta^{\prime} ) \geq
1}\frac{%
\Vert\Sigma( \theta) -\Sigma( \theta
^{\prime
} ) \Vert^{2}}{H ( \theta,\theta^{\prime} )
}\geq
cka^{2}
\]
for some constant $c>0$.
\end{pf}

\section{Simulation study}
\label{simulation.sec}

We now turn to the numerical performance of the proposed tapering estimator
and compare it with that of the banding estimator of Bickel and Levina
(\citeyear{BL08a}). In the numerical study, we shall consider estimating a covariance
matrix in the parameter space $\mathcal{F}_{\alpha}$ defined in (\ref
{paraspace}). Specifically, we consider the covariance matrix $\Sigma
= ( \sigma_{ij} ) _{1\leq i,j\leq p}$ of the form
%
%
\begin{equation}
\sigma_{ij} = \cases{
1, &\quad $1\leq i=j \leq p$, \cr
\rho\vert i-j \vert^{-(\alpha+1)}, &\quad $1\leq i\neq j\leq p$.}
\end{equation}
Note that this is a Toeplitz matrix. But we do not assume that the
structure is known and do not use the information in any estimation
procedure.

The banding estimator in (\ref{bandestimator}) depends on the choice
of $k$.
An optimal tradeoff of $k$ is $k\asymp( n/\log p )
^{1/ (
2\alpha+2 ) }$ as discussed in Section \ref{sec.discussion}. See
Bickel and Levina (\citeyear{BL08a}). The tapering estimator (\ref{estimator}) also
depends on $k$ for which the optimal tradeoff is $k\asymp n^{1/ (
2\alpha+1 ) }$. In our simulation study, we choose
%
\begin{table}
\tabcolsep=0pt
\caption{The average errors under the spectral norm of the banding estimator
(BL) and the tapering estimator (CZZ) over $100$ replications.
The cases
where the tapering estimator underperforms the banding estimator are
highlighted in italic}
\label{table}
\begin{tabular*}{\tablewidth}{@{\extracolsep{\fill}}ld{4.0}cccccccccc@{}}
\hline
& & \multicolumn{2}{c}{$\bolds{\alpha=0.1}$} &
\multicolumn{2}{c}{$\bolds{\alpha=0.2}$} &
\multicolumn{2}{c}{$\bolds{\alpha=0.3}$} & \multicolumn{2}{c}{$\bolds{\alpha=0.4}$}
& \multicolumn{2}{c@{}}{$\bolds{\alpha=0.5}$}\\[-4pt]
& & \multicolumn{2}{c}{\hrulefill} &
\multicolumn{2}{c}{\hrulefill} &
\multicolumn{2}{c}{\hrulefill} & \multicolumn{2}{c}{\hrulefill} &
\multicolumn{2}{c@{}}{\hrulefill}\\
$\bolds p$ & \multicolumn{1}{c}{$\bolds n$} & \textbf{BL} & \textbf{CZZ} & \textbf{BL} & \textbf{CZZ}
& \textbf{BL} & \textbf{CZZ} & \textbf{BL} & \textbf{CZZ} & \textbf{BL} & \textbf{CZZ}\\
\hline
\phantom{0}250 & 250
& 2.781 & 2.706 & 2.291 & 2.023 & 1.762 & 1.684 & 1.618 & 1.517 &
\textit{1.325} & \textit{1.507}\\
& 500
& 2.409 & 2.302 & 1.898 & 1.575 & 1.562 & 1.204 & 1.361 & 1.185 & 1.080
& 0.822\\
& 1000
& 2.029 & 1.685 & 1.631 & 1.361 & 1.289 & 1.018 & 1.056 & 0.795 & 0.911
& 0.859\\
& 2000
& 1.706 & 1.153 & 1.369 & 1.122 & 1.106 & 0.908 & 0.878 & 0.655 & 0.715
& 0.542\\
& 3000 & 1.522 & 0.926 & 1.242 & 0.896 & 0.983 & 0.798 & 0.810 & 0.658
& 0.645 & 0.482\\
[4pt]
\phantom{0}500 & 250 & 3.277 & 2.914 & 2.609 & 2.097 & 1.961 & 1.788 & 1.745 & 1.610 &
\textit{1.392} & \textit{1.571}\\
& 500
& 2.901 & 2.598 & 2.199 & 1.683 & 1.751 & 1.256 & 1.475 & 1.234 & 1.152
& 0.865\\
& 1000
& 2.539 & 2.197 & 1.942 & 1.472 & 1.481 & 1.064 & 1.178 & 0.843 & 0.984
& 0.917\\
& 2000 & 2.263 & 1.726 & 1.669 & 1.326 & 1.293 & 0.965 & 1.067 & 0.700
& 0.866 & 0.569\\
& 3000 & 2.066 & 1.379 & 1.538 & 1.154 & 1.220 & 0.874 & 0.919 & 0.696
& 0.781 & 0.503\\
[4pt]
1000 & 250 & 3.747 & 3.086 & 2.873 & 2.223 & 2.385 & 1.842 & 1.833 & 1.694 &
\textit{1.449} & \textit{1.643}\\
& 500 & 3.370 & 2.735 & 2.635 & 1.768 & 1.906 & 1.334 & 1.565 & 1.297 &
1.203 & 0.925\\
& 1000 & 3.097 & 2.437 & 2.315 & 1.536 & 1.741 & 1.121 & 1.382 &
0.883 & 1.037 & 0.936\\
& 2000 & 2.730 & 2.177 & 2.011 & 1.392 & 1.523 & 1.006 & 1.156 & 0.722
& 0.920 & 0.591\\
& 3000 & 2.589 & 1.968 & 1.865 & 1.264 & 1.374 & 0.911 & 1.072 & 0.723
& 0.834 & 0.523\\
[4pt]
2000 & 250 & 4.438 & 3.177 & 3.107 & 2.300 & 2.511 & 1.956 & 1.903 & 1.744 &
\textit{1.484} & \textit{1.736}\\
& 500 & 3.969 & 2.800 & 2.868 & 1.841 & 2.030 & 1.383 & 1.638 & 1.356 &
1.239 & 0.940\\
& 1000 & 3.538 & 2.531 & 2.551 & 1.599 & 1.866 & 1.158 & 1.452 &
0.912 & 1.074 & 0.973\\
& 2000 & 3.242 & 2.353 & 2.248 & 1.434 & 1.649 & 1.031 & 1.224 & 0.751
& 0.955 & 0.611\\
& 3000 & 3.025 & 2.219 & 2.101 & 1.302 & 1.566 & 0.929 & 1.141 & 0.743
& 0.868 & 0.541\\
[4pt]
3000 & 250
& 4.679 & 3.219 & 3.230 & 2.358 & 2.576 & 1.995 & 1.931 & 1.797 &
\textit{1.494} & \textit{1.776}\\
& 500
& 4.214 & 2.887 & 2.991 & 1.890 & 2.282 & 1.419 & 1.664 & 1.384 & 1.463
& 0.971\\
& 1000
& 3.901 & 2.575 & 2.674 & 1.633 & 1.933 & 1.186 & 1.482 & 0.929 & 1.224
& 0.990\\
& 2000
& 3.488 & 2.395 & 2.452 & 1.451 & 1.717 & 1.049 & 1.254 & 0.768 & 0.965
& 0.619\\
& 3000
& 3.336 & 2.278 & 2.288 & 1.321 & 1.632 & 0.948 & 1.172 & 0.750 & 0.880
& 0.549\\
\hline
\end{tabular*}
\end{table}
$k=
\lfloor
( n/\log p ) ^{1/ ( 2\alpha+2 ) } \rfloor
$ for the
banding estimator and $k= \lfloor n^{1/ ( 2\alpha+1 )
} \rfloor$ for the tapering estimator.

A range of parameter values for $\alpha$, $n$ and $p$ are considered.
Specifically, $\alpha$ ranges from $0.1$ to $0.5$, the sample size $n$
ranges from $250$ to $3000$ and the dimension $p$ goes from $250$ to $3000$.
We choose the value of $\rho$ to be $\rho=0.6$ so that all matrices are
nonnegative definite and their smallest eigenvalues are close to $0$.
Table \ref{table} reports the average errors under the spectral norm over $100$
replications for the two procedures. The cases where the tapering estimator
underperforms the banding estimator are highlighted in boldface. Figure
\ref{risk-comp} plots the ratios of the average errors of the banding estimator
to the corresponding average errors of the tapering estimator for
$\alpha
=0.1,0.2,0.3$ and $0.5$. The case of $\alpha=0.4$ is similar to the
case of
$\alpha=0.3$.

It can be seen from Table \ref{table} and Figure \ref{risk-comp} that the
tapering estimator outperforms the banding estimator in 121 out of 125
cases. For the given dimension $p$, the ratio of the average error of the
banding estimator to the corresponding average error of the tapering
estimator tends to increase as the sample size $n$ increases. The tapering
estimator fails to outperform the banding estimator only when $\alpha=0.5$
and $n=250$ in which case the values of $k$ are small for both estimators.
\begin{remark}
We have also carried out additional simulations for larger values
of $\alpha$ with the same sample sizes and dimensions. The
performance of the tapering and abnding estimators
are similar. This is mainly dur to the fact that the values of $k$ for
both estimators are very small for large $\alpha$ when $n$ and $p$ are
only moderately large.
\end{remark}

\section{Proofs of auxiliary lemmas}
\label{sec.proofs}

In this section we give proofs of auxiliary lemmas stated and used in
Sections \ref{operatornorm.sec}--\ref{inverse.sec}.
\begin{pf*}{Proof of Lemma \protect\ref{est}}
Without loss of generality we assume that $i\leq j$. The set $ \{
i,j \} $ is contained in the set $ \{ l,\ldots
,l+k_{h}-1 \} $
if and only if $l\leq i\leq j\leq l+k_{h}-1$, that is, $j-k_{h}+1\leq
l\leq i$.
Note that $\operatorname{Card} \{ l\dvtx j-k_{h}+1\leq l\leq i
\} = (
i- ( j-k_{h}+1 ) +1 ) _{+}=(k_{h}-|i-j|)_{+}$, then
$\operatorname{Card} \{ l\dvtx \{ i,j \} \subset\{
l,\ldots
,l+k_{h}-1 \} \} =(k_{h}-|i-j|)_{+}$. Similarly, we have
%
%
\begin{figure}

\includegraphics{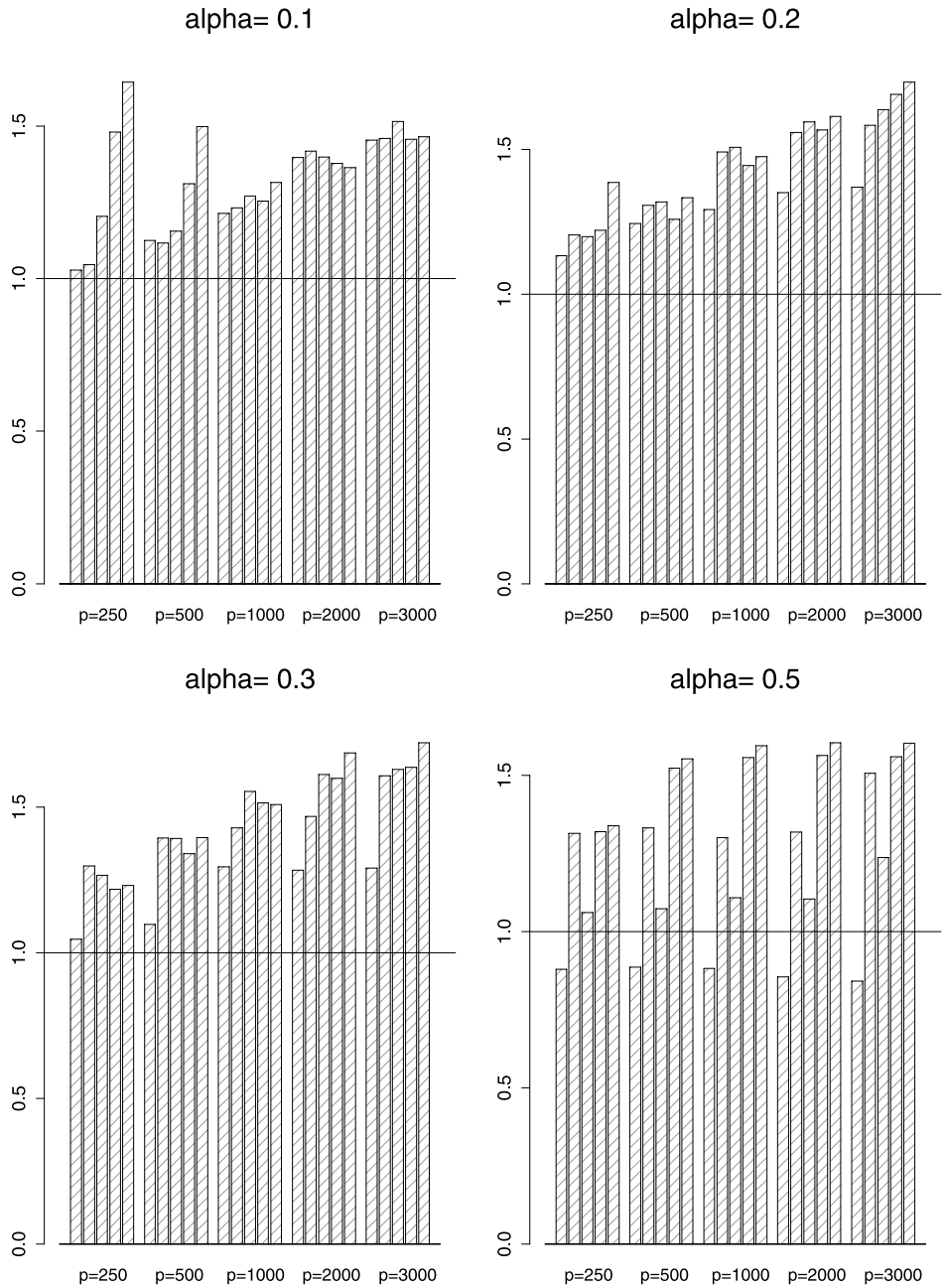}

\caption{The vertical bars represent the ratios of the
average error of the banding estimator to the corresponding average
error of
the tapering estimator. The higher the bar the better the relative
performance of the tapering estimator. For each value of $p$ the bars are
ordered from left to right by the sample sizes $(n=250$ to $3000)$.}
\label{risk-comp}
\end{figure}
$\operatorname{Card} \{ l\dvtx \{ i,j \} \subset\{
l,\ldots,l+k-1 \}
\} =(k-|i-j|)_{+}$. Thus we have
\begin{eqnarray*}
kw_{ij} &=&(k-|i-j|)_{+}-(k_{h}-|i-j|)_{+} \\
&=&\operatorname{Card} \bigl\{ l\dvtx \{ i,j \} \subset
\{ l,\ldots
,l+k-1 \} \bigr\} \\
&&{}-\operatorname{Card} \bigl\{ l\dvtx \{
i,j \} \subset
\{ l,\ldots,l+k_{h}-1 \} \bigr\}.
\end{eqnarray*}
\upqed\end{pf*}
\begin{pf*}{Proof of Lemma \protect\ref{estbias}}
Without loss of generality we assume that $p$ is divisible by $m$. Recall
that $M_{l}^{ ( m ) }= ( \tilde{\sigma}_{ij}I \{
l\leq
i<l+m,l\leq j<l+m \} ) _{p\times p}$. Note that
$M_{l}^{ (
m ) }$ is empty when $l\leq1-m$, and has at least one nonzero entry
when $l\geq2-m$. Set $\delta_{l}^{(m)}=M_{l}^{(m)}-\mathbb{E}M_{l}^{(m)}$
and $S^{ ( m ) }=\sum_{l=2-m}^{p}M_{l}^{ ( m )
}$. It
follows from (\ref{estimator}) that
%
%
\begin{equation} \label{pdivided}
\bigl\Vert S^{(m)}-\mathbb{E}S^{(m)} \bigr\Vert\leq
\sum_{l=1}^{m} \biggl\Vert\sum_{-1\leq j < p/m}\delta
_{jm+l}^{(m)} \biggr\Vert.
\end{equation}
Since $\delta_{jm+l}^{(m)}$ are disjoint diagonal blocks over $-1\leq
j<p/m$%
, we have
%
%
\begin{eqnarray}\label{bound22}
\bigl\Vert S^{(m)}-\mathbb{E}S^{(m)} \bigr\Vert &\leq& m\max_{1\leq
l\leq
m} \biggl\Vert\sum_{-1\leq j < p/m}\delta_{jm+l}^{(m)} \biggr\Vert
\nonumber\\[-8pt]\\[-8pt]
&\leq&
m\max_{1-m\leq l\leq p} \bigl\Vert\delta_{l}^{(m)} \bigr\Vert.\nonumber
\end{eqnarray}
Since $\delta_{l}^{(k_{h})}$ and $\delta_{l}^{(k)}$ are all
sub-blocks of
certain matrix $\delta_{l}^{(k)}$ with $1\leq l\leq p-k+1$, Lemma \ref
{estbias} now follows immediately from equations (\ref{bound22}) and
(\ref%
{estimator}).
\end{pf*}
\begin{pf*}{Proof of Lemma \protect\ref{estbiasbd}}
For any $m\times m$ symmetric matrix $A$, we have
\begin{eqnarray*}
\vert u^{T}Au \vert- \vert v^{T}Av \vert
&\leq&
\vert
u^{T}Au-v^{T}Av \vert= \vert( u-v ) ^{T}A (
u+v ) \vert\\
&\leq&\Vert u-v \Vert\Vert
A \Vert\Vert u+v \Vert.
\end{eqnarray*}
Let $S_{1/2}^{m-1}$ be a $1/2$ net of the unit sphere $S^{m-1}$ in the
Euclidean distance in ${R}^{m}$. We have%
\begin{eqnarray*}
\Vert A \Vert&\leq&{\sup_{u\in S^{m-1}}} \vert
u^{T}Au \vert\leq{\sup_{u\in S_{1/2}^{m-1}}} \vert
u^{T}Au \vert+\frac{1}{2} \Vert A \Vert\frac{3}{2}\\
&=&{\sup_{u\in S_{1/2}^{m-1}}} \vert u^{T}Au \vert+\frac{3}{4}
\Vert A \Vert,
\end{eqnarray*}
which implies\vspace*{-2pt} $ \Vert A \Vert\leq{4\sup_{u\in
S_{1/2}^{m-1}}} \vert u^{T}Au \vert$. Since we are allowed
to pack\break
$\operatorname{Card} ( S_{1/2}^{m-1} ) $ balls\vspace*{1pt} of radius
$1/4$ into a $%
1+1/4$ ball in ${R}^{m}$, volume comparison yields
\[
(1/4)^{m}\operatorname{Card} ( S_{1/2}^{m-1} ) \leq(5/4)^{m},
\]
that is, $\operatorname{Card} ( S_{1/2}^{m-1} ) \leq
5^{m}$. Thus there
exist $\mathbf{v}_{1},\mathbf{v}_{2},\ldots,\mathbf{v}_{5^{m}}\in S^{m-1}$
such that%
\[
\Vert A \Vert\leq{4\sup_{j\leq5^{m}}} \vert
v_{j}^{T}Av_{j} \vert\qquad\mbox{for all }m\times m\mbox{ symmetric }A.
\]
This one-step approximation argument is similar to the proof of Proposition
4.2(ii) in Zhang and Huang (\citeyear{ZH08}).

Let $\mathbf{X}_{1},\ldots,\mathbf{X}_{n}$ be i.i.d. $p$-vectors
with $%
\mathbb{E} ( \mathbf{X}_{1}\mathbf{-\mu} ) (
\mathbf{X}_{1}%
\mathbf{-\mu} ) ^{T}=\Sigma$. Under the sub-Gaussian assumption
in (%
\ref{subGau}) there exists $\rho>0$ such that
\[
\mathbb{P} \{ \mathbf{v}^{T}(\mathbf{X}_{i}-\mathbb{E}\mathbf
{X}_{i})(%
\mathbf{X}_{i}-\mathbb{E}\mathbf{X}_{i})^{T}\mathbf{v}>x \}
\leq
e^{-x\rho/2}\qquad\mbox{for all }x>0\mbox{ and }\Vert\mathbf{v}\Vert=1,
\]
which implies $\mathbb{E} ( t\mathbf{v}^{T}(\mathbf
{X}_{i}-\mathbb{E}%
\mathbf{X}_{i})(\mathbf{X}_{i}-\mathbb{E}\mathbf{X}_{i})^{T}\mathbf
{v}%
) <\infty$ for all $t<\rho/2$ and $\Vert\mathbf{v}\Vert=1$, then
there exists $\rho_{1}>0$ such that
\[
\mathbb{P} \Biggl\{ \Biggl\vert\frac{1}{n}\sum_{i=1}^{n}\mathbf
{v}^{T} [ (%
\mathbf{X}_{i}-\mathbb{E}\mathbf{X}_{i})(\mathbf{X}_{i}-\mathbb
{E}\mathbf{X}%
_{i})^{T}-\Sigma] \mathbf{v} \Biggr\vert>x \Biggr\} \leq
e^{-nx^{2}\rho_{1}/2}
\]
for all $0<x<\rho_{1}$ and $\Vert\mathbf{v}\Vert=1$. [See, e.g., Chapter
2 in Saulis and Statulevi\v{c}ius (\citeyear{SS91}).] Thus we have
\begin{eqnarray*}
&&\mathbb{P} \Bigl\{ \max_{1\leq l\leq p-m+1} \bigl\Vert
M_{l}^{(m)}-\mathbb{E}%
M_{l}^{(m)} \bigr\Vert>x \Bigr\} \\
&&\qquad\leq\sum_{1\leq l\leq
p-m+1}\mathbb{P}%
\bigl\{ \bigl\Vert M_{l}^{(m)}-\mathbb{E}M_{l}^{(m)} \bigr\Vert
>x \bigr\} \\
&&\qquad\leq 2p5^{m}\sup_{\mathbf{v}_{j},l}\mathbb{P}\bigl\{\bigl|\mathbf{v}%
_{j}^{T}\bigl(M_{l}^{(m)}-\mathbb{E}M_{l}^{(m)}\bigr)\mathbf{v}_{j}\bigr|>x\bigr\} \\
&&\qquad\leq 2p5^{m}\exp( -nx^{2}\rho_{1}/2 ). 
\end{eqnarray*}
\upqed\end{pf*}
\begin{pf*}{Proof of Lemma \protect\ref{dffbd}}
Set $v= ( 1 \{ k_{h}\leq i\leq k \} ) $ and let
\[
( w_{i} ) = [ \Sigma( \theta) -\Sigma
(
\theta^{\prime} ) ] v.
\]
Note that there are exactly $H ( \theta,\theta^{\prime} ) $
number of $w_{i}$ such that $ \vert w_{i} \vert=\tau
k_{h}a$, and $%
\Vert v \Vert_{2}^{2}=k_{h}$. This implies%
\begin{eqnarray*}
\Vert\Sigma( \theta) -\Sigma( \theta
^{\prime
} ) \Vert^{2} &\geq& \frac{ \Vert[ \Sigma
( \theta
) -\Sigma( \theta^{\prime} ) ] v
\Vert_{2}^{2}%
}{ \Vert v \Vert_{2}^{2}}
\geq\frac{H(\theta,\theta
^{\prime
})\cdot( \tau ka ) ^{2}}{k_{h}}\\
&=& H(\theta,\theta^{\prime
})\cdot
\tau^{2}k_{h}a^{2}.
\end{eqnarray*}
\upqed\end{pf*}
\begin{pf*}{Proof of Lemma \protect\ref{affbd}}
When $H ( \theta,\theta^{\prime} ) =1$, we will show
\begin{eqnarray*}
\Vert P_{\theta^{\prime}}-P_{\theta} \Vert_{1}^{2}
&\leq&
2K ( P_{\theta^{\prime}}|P_{\theta} ) \\
&=& 2n \biggl[ \frac
{1}{2}%
\operatorname{tr} ( \Sigma( \theta^{\prime} ) \Sigma^{-1} (
\theta
) ) -\frac{1}{2}\log\det( \Sigma( \theta
^{\prime
} ) \Sigma^{-1} ( \theta) ) -\frac
{p}{2} \biggr] \\
&\leq& n\cdot cka^{2}
\end{eqnarray*}
for some small $c>0$, where $K ( \cdot|\cdot) $ is the
Kullback--Leibler divergence and the first inequality follows from the
well-known Pinsker's inequality [see, e.g., Csisz\'{a}r (\citeyear{C67})]. This
immediately
implies the $L_{1}$ distance between two measures is bounded away from $1$,
and then the lemma follows. Write
\[
\Sigma( \theta^{\prime} ) =D_{1}+\Sigma( \theta
).
\]
Then%
\[
\frac{1}{2}\operatorname{tr} ( \Sigma( \theta^{\prime} ) \Sigma
^{-1} ( \theta) ) -\frac{p}{2}=\frac{1}{2}\operatorname{tr} (
D_{1}\Sigma^{-1} ( \theta) ).
\]
Let $\lambda_{i}$ be the eigenvalues of $D_{1}\Sigma^{-1} (
\theta
) $. Since $D_{1}\Sigma^{-1} ( \theta) $ is
similar to the
symmetric matrix $\Sigma^{-1/2} ( \theta) D_{1}\Sigma
^{-1/2} ( \theta) $, and%
\begin{eqnarray*}
\Vert\Sigma^{-1/2} ( \theta) D_{1}\Sigma
^{-1/2} (
\theta) \Vert
&\leq& \Vert\Sigma^{-1/2}(\theta
) \Vert\Vert D_{1} \Vert\Vert\Sigma
^{-1/2} (
\theta) \Vert\\
&\leq& c_{1} \Vert D_{1} \Vert
\leq
c_{1} \Vert D_{1} \Vert_{1}\leq c_{2}ka,
\end{eqnarray*}
then all eigenvalues $\lambda_{i}$'s are real and in the interval
$ [
-c_{2}ka,c_{2}ka ] $, where $ka=k\cdot k^{- ( \alpha
+1 )
}=k^{-\alpha}\rightarrow0$. Note that the Taylor expansion yields%
\[
\log\det( \Sigma( \theta^{\prime} ) \Sigma
^{-1} (
\theta) ) =\log\det\bigl( I+D_{1}\Sigma^{-1} (
\theta
) \bigr) =\operatorname{tr} ( D_{1}\Sigma^{-1} ( \theta)
)-R_{3},
\]
where
\[
R_{3}\leq c_{3}\sum_{i=1}^{p}\lambda_{i}^{2}\qquad\mbox{for some }c_{3}>0.
\]
Write $\Sigma^{-1/2} ( \theta) =UV^{1/2}U^{T}$, where $UU^{T}=I$
and $V$ is a diagonal matrix. It follows from the fact that the Frobenius
norm of a matrix remains the same after an orthogonal transformation that
\begin{eqnarray*}
\sum_{i=1}^{p}\lambda_{i}^{2} &=& \Vert\Sigma^{-1/2} (
\theta) D_{1}\Sigma^{-1/2} ( \theta) \Vert
_{F}^{2}
\leq\Vert V \Vert^{2}\cdot\Vert
U^{T}D_{1}U \Vert_{F}^{2}\\
&=& \Vert\Sigma^{-1} (\theta)
\Vert^{2}\cdot\Vert D_{1} \Vert_{F}^{2}\leq
c_{4}ka^{2}.
\end{eqnarray*}
\upqed\end{pf*}

\section*{Acknowledgments}
The authors would like to thank
James X. Hu for assistance in carrying out the simulation study in
Section %
\ref{simulation.sec}. We also thank the Associate Editor and three referees
for thorough and useful comments which have helped to improve the
presentation of the paper.

\printaddresses

\end{document}